\documentclass[10pt,a4paper]{article}
\usepackage[utf8]{inputenc}
\usepackage{amsmath}
\usepackage{amsfonts}
\usepackage{amssymb}
\usepackage{makeidx}
\usepackage{graphicx}
\graphicspath{ {./} }
\usepackage{listings}
\usepackage{tocloft}
\usepackage{verbatim}
\usepackage{xcolor}
\usepackage{calc}
\usepackage{CJKutf8}
\usepackage[hidelinks]{hyperref}

\usepackage{epigraph}
\newtheorem{le[mm]a}[theorem]{Le[mm]a}

\newcommand{\qed}{\nobreak \ifvmode \relax \else
      \ifdim\lastskip<1.5em \hskip-\lastskip
      \hskip1.5em plus0em minus0.5em \fi \nobreak
      \vrule height0.75em width0.5em depth0.25em\fi}
      
       \usepackage[polutonikogreek,english]{babel}
    
\usepackage[normalem]{ulem}

\usepackage{tikz, tkz-tab, pgfkeys}
\usetikzlibrary{calc,decorations.pathmorphing,shapes}

  \newcounter{exer}[section]

    \makeindex

 \title{Johann Heinrich Lambert's memoir \\
\emph{Theorie der Parallellinien}: \\
A review with commentary}
 \author{Athanase Papadopoulos and Guillaume Théret}
 
\begin{document}
\maketitle

\epigraph{\itshape 
 Yet there have been and still are geometers and philosophers, and even some of the most distinguished, who doubt whether the whole universe, or to speak more widely, the whole of being, was only created in Euclid's geometry; they even dare to dream that two parallel lines, which according to Euclid can never meet on earth, may meet somewhere in infinity. I have come to the conclusion that, since I can't understand even that, I can't expect to understand about God. I acknowledge humbly that I have no faculty for settling such questions, I have a Euclidean earthly mind, and how could I solve problems that are not of this world? And I advise you never to think about it either, my dear Alyosha, especially about God, whether He exists or not. All such questions are utterly inappropriate for a mind created with an idea of only three dimensions. (Fyodor Dostoyevsky, \emph{The Brothers Karamazov}, 
transl.  from the Russian
by Constance Garnett,
The Lowell Press,
New York,
p. 257 [Ivan Karamazov speaking].)}

\vspace*{\fill} \epigraph{\itshape 
 Like the parallel postulate, the
postulate that mathematics will survive has been stripped of its ``evidence";
but, while the former is no longer necessary, we would not be able to get on
without the latter. (Andr\'e Weil, \emph{The Future of Mathematics}, 1946.)}

\vfill\eject

 \noindent{\bf Abstract} We review the memoir \emph{Theorie der Parallellinien} by Johann Heinrich Lambert, written in 1766. Lambert, a victim of the prejudices of his time, conceived this memoir as an attempt to prove the so-called parallel postulate of Euclid's \emph{Elements}, and consequently, the non-existence of the geometry that we now call hyperbolic geometry. In fact, by developing the foundations of a geometry obtained by replacing the\index{parallel postulate} parallel postulate with its negation while keeping Euclid's other postulates unchanged, Lambert was hoping to arrive at a contradiction. Of course, he failed in his endeavor, but these attempts at proving the parallel postulate implicitly contain, without Lambert having foreseen it, fundamental results of hyperbolic geometry, the discovery of which, by Lobachevsky, Bolyai and Gauss, was not to take place until the following century. Thus, Lambert's memoir (which he did not intend to publish but which was eventually published in 1895) constitutes one of the founding texts of non-Euclidean geometry. 

      Spherical geometry is one of the three geometries of constant curvature, the other two being Euclidean geometry and hyperbolic geometry. In this sense, along with hyperbolic geometry, spherical geometry constitutes one of the two non-Euclidean geometries. In fact, Lambert, like Lobachevsky and others after him, understood the deep relationships between the three geometries: Euclidean, spherical, and hyperbolic, in particular the formal and the more profound analogies between the trigonometric formulae, the properties of birectangular isosceles quadrilaterals and of trirectangular quadrilaterals, the monotonicity properties (which can be formulated in terms of convexity properties) which hold in opposite senses in spherical and hyperbolic geometry which at some points he calls a sphere of imaginary radius.
It is for these reasons that we decided to include in this volume, dedicated to spherical geometry, a chapter on this important memoir by Lambert, trying to highlight its most important ideas.

This paper will appear as a chapter in the book \emph{Spherical Geometry in the Eighteenth Century I: Euler, Lagrange and Lambert}, ed. R. Caddeo and A. Papadopoulos, Springer Nature Switzerland.
 
 \medskip


  
 \medskip

  \section{Introduction: The theory of parallels}
  
  The expression ``theory of parallel lines"\index{parallel postulate} in the title of Lambert's memoir \emph{Theorie der Parallellinien}, published in \cite{Engel-Staeckel} (see \cite{2012-Lambert} for a French translation), refers to the various attempts that were made, over a period that lasted almost two and a half millennia, to prove the so-called parallel postulate\index{fifth postulate} of Euclid's\index{Euclid} \emph{Elements}.
We recall\index{Euclid!\emph{Elements}} that this postulate is the statement placed  in Book 1 of J. L. Heiberg's classic edition of the \emph{Elementa} as postulate 5  \cite[p. 8]{Heiberg} and which, in Heath's English edition, based on Heiberg's, says the following  \cite[Vol. I, p. 202]{Heath3}: 
\begin{quote}\small
That, if a straight line falling on two straight lines makes the interior angles
on the same side less than two right angles, the two straight lines, if produced
indefinitely, meet on that side on which are the angles less than the two right
angles.
\end{quote}

 In some other versions of the \emph{Elements}, this statement is included as the eleventh axiom. Lambert uses the term ``eleventh axiom." 
 
  Euclid\index{Euclid} makes a distinction between a postulate\index{postulate} and an axiom.\index{axiom} This difference is discussed\index{parallel postulate} in detail by Proclus\footnote{\label{f:Proclus} Proclus (ca. 412-485 AD) was a Neoplatonist philosopher who studied in Alexandria. His work is an invaluable source of knowledge\index{Proclus of Lycia} for historians of Greek mathematics. We have included some information on him in \S \ref{excursus:parallel-lines} of the present chapter.}  in his \emph{Commentaries on Book I of Euclid's Elements} \cite{Proclus}.  To put it in simple words,  a postulate asserts that a certain geometric construction is possible, while an axiom may concern properties (such as magnitudes) that are not geometric. The distinction between ``postulate" and ``axiom" has disappeared in modern axiomatics, that is, those of Peano, Veronese, Pasch, Hilbert, etc., and we shall not make it here. Thus, in this chapter, we will use the words ``postulate" or ``axiom" interchangeably, and therefore also the expressions  ``parallel axiom", ``fifth axiom", ``parallel postulate", ``eleventh postulate".
   
The history of this ``theory of parallels" covers a period ranging from Greek antiquity (starting before Plato) until the time (around 1868) when the works of Lobachevsky\index{Lobachevsky, Nikolai Ivanovich}, Bolyai\index{Bolyai, J\'anos} and Gauss\index{Gauss, Carl Friedrich} on non-Euclidean\index{non-Euclidean geometry} geometry started to be considered by the mathematical community as valid. These discoveries showed that Euclid's parallel postulate cannot be deduced from the other postulates, and it gradually put an end to the very numerous attempts to prove this postulate.

Mathematicians addressed the theory of parallels from several facets. First, there are geometers who tried to deduce the fifth postulate from Euclid's other postulates without introducing a new one. We know a posteriori that these attempts were doomed to failure. Other authors modified the definition\index{parallel lines!definition!Euclid} of parallelism,\footnote{In Euclid's \emph{Elements}, this definition is given in Book I (Definition 25): \emph{Parallel straight lines are straight lines which,
being in the same plane and being produced indefinitely in
both directions, do not meet one another in either direction}. \cite[Vol. I, p. 154]{Heath3}.} believing that this would avoid the difficulty, and these also, in general, achieved nothing. But there are others who, finding the postulate complicated and less intuitive than the other postulates, tried to remedy this by replacing it with a simpler and more intuitive one: equidistance of parallel lines, existence of translations, existence of homothetic triangles, etc. Some of them succeeded in their endeavor and obtained valuable results. In Section \ref{excursus:parallel-lines} of the present chapter, we review several works that have been done on this topic by mathematicians over a period from Greek Antiquity to the 19th century, passing through the works of the Arab mathematicians of the Middle Ages.

Lambert holds a special place in this story. His work, although it predates the discovery of non-Euclidean geometry,\index{non-Euclidean geometry} implicitly contains results from that theory, and this is what we shall try to explain in the following pages.

It is conceivable that the question of parallels is, among all the questions that have been addressed to date by prominent mathematicians, the one about which the volume of existing writings is the largest. Thousands of articles and hundreds of books, of various styles and scopes and spread over a period of more than two thousand years, are devoted to this theory. We may quote here Henri Poincar\'e who,\index{Poincar\'e, Henri} in an exposition of the question published in 1891 \cite{Poincare1891}, writes: ``The amount of effort that has been expended on this chimerical hope is truly unthinkable.''
An article by the French mathematician Georges Brunel\index{Brunel, Grorges} (1854-1900) published in 1888 in the \emph{M\'emoires de la Soci\'et\'e de Sciences Physiques et Naturelles de Bordeaux} \cite{Brunel} mentions one hundred and thirty-four attempts to prove this postulate. The book by Paul St\"ackel\index{St\"ackel, Paul} and Friedrich Engel\index{Engel, Friedrich} titled \emph{Die Theorie der Parallellinien von Euklid bis auf Gauss, eine Urkundensammlung zur Vorgeschichte der nicht-euklidischen Geometrie} (Theory of parallel lines from Euclid to Gauss: A collection of sources for the prehistory of non-Euclidean geometry) \cite{Engel-Staeckel}, published in 1895, contains, besides the text by Lambert which is the subject of the present chapter, excerpts of works by Euclid,\index{Euclid} Wallis,\index{Wallis, John} Saccheri\index{Saccheri, Giovanni Girolamo}, Gauss\index{Gauss, Carl Friedrich}, Schweikart\index{Schweikart, Ferdinand Karl} and Taurinus\index{Taurinus, Franz Adolph} on the theory of parallels. The bibliography of this book includes three hundred and twenty items about the theory of parallels, published in the period 1482-1837, a bibliography which remains minimalist compared to all that has been written on this subject.  St\"ackel and Engel, for example, left aside the very important contribution of the Greeks to the question, as well as the no less important contribution of the Arabs, the latter having been, it is true, very little known at the time -- with a few notable exceptions; for instance, the work of Na\d{s}\=\i r al-D\=\i n al-\d{T}\=us\=\i \ (13th century)\index{Tusi@al-\d{T}\=us\=\i , Na\d{s}\=\i r al-D\=\i n} was known and quoted by Wallis and Saccheri. The book by J.-C. Pont, \cite{Pont}, \emph{L'aventure des parall\`eles} (The Adventure of Parallels), published in 1986, presents a very good overview of this subject, including a section on the work of the Greeks and another one on that of the Arabs.

The plan of the rest of this chapter is the following:
Section \ref{s:structure} contains an outline of the structure of Lambert's memoir. It is followed by an analysis of the three parts of this memoir, made in \S \ref{s:Lambert-Part-1}, \ref{s:Lambert-Part-2} and \ref{s:Lambert-Part-3}. After this, we have included two excursuses. The first one, in \S  \ref{excursus:parallel-lines}, contains some milestones on the history of parallel lines, from
Greek antiquity until the modern period, and the second one, in \S \ref{excursus:quadrilaterals}, is concerned with two sorts of quadrilaterals that play a central role in the theory of parallels, namely, the trirectangular  and birectangular isosceles quadrilaterals.

\section{The structure of Lambert's memoir}\label{s:structure}

Lambert's memoir, \emph{Theorie der Parallellinien}, is divided into three parts:

\begin{enumerate}

 \item \S 1 to 11: A description of the problem of parallels and a discussion of Euclid's\index{Euclid} fifth postulate and its place among the propositions and other axioms of the \emph{Elements}.\index{Euclid!\emph{Elements}} Lambert discusses there the difficulties posed by this postulate and recalls commentaries and attempts by several of his predecessors to prove it.
 
  \item \S 12 to 26:  
     Some propositions of neutral geometry, that is, of Euclidean\index{Euclidean geometry} geometry deprived of the parallel axiom. Neutral geometry is defined axiomatically by taking the Euclidean system of axioms  and postulates with the exception of the parallel postulate. Lambert hopes that these propositions of neutral geometry may  be used to prove the fifth postulate.
 
   \item \S 27 to 88:  Lambert's attempts to prove the parallel postulate, starting from the principle that a negation of the parallel axiom leads to an absurdity. It is in this third part that we find a good number of propositions that are part of the field that was later called hyperbolic geometry.

   We shall survey each part of Lambert's memoir  separately.
   \end{enumerate}

\section{Part 1 of Lambert's memoir}\label{s:Lambert-Part-1}

\begin{figure}
\centering
\includegraphics[width=0.8\linewidth]{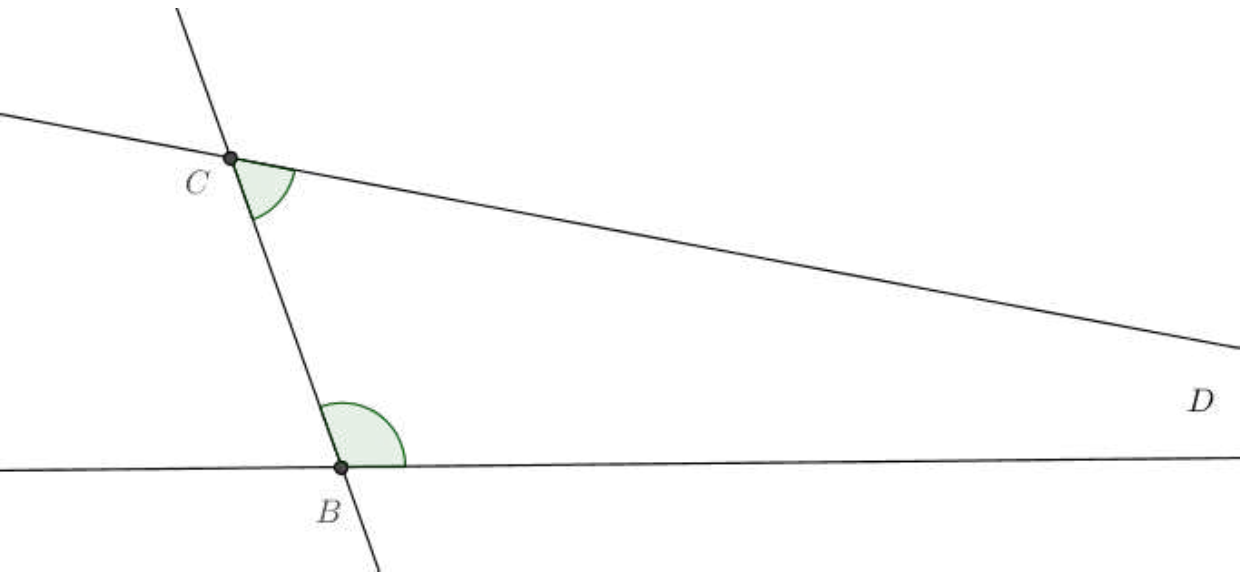}
\caption{\small The figure used in Lambert's statement of Euclid's parallel axiom} \label{fig:2-Lambert1}
\end{figure}

In this part, Lambert undertakes an exposition of what he calls ``the general problem of parallel lines", a problem which, he says, constitutes ``a difficulty that appeared at the dawn of geometry and which, since the time of Euclid,\index{Euclid} has struck all those who refused to blindly believe, following others, the theories of this science and who on the contrary wanted to be convinced by reason and never to depart from the rigor that they found in most of the proofs."
 The difficulty, which, he says, jumps out at every reader\index{Euclid!\emph{Elements}} of the \emph{Elements}, lies in the eleventh axiom. He states this axiom in the following form (Figure \ref{fig:2-Lambert1}):\footnote{In the following sentence, and especially in the figure, the letter $D$ indicates a direction rather than a point.}

   \begin{quote}\small
   
 When two lines $CD, BD$ are cut by a third one $BC$ and when the two internal angles $\widehat{DCB}$ and $\widehat{DBC}$, taken together, are smaller than two right angles, then the two lines $CD, BD$ intersect on the side of $D$, 
that is, on the side where we find these angles.\footnote{Compare with the statement in Euclid's \emph{Elements}, Heath's translation \cite[Vol. I, p. 202]{Heath3} which we recalled in the Introduction to the present chapter.}
 
 \end{quote}

 Lambert declares that when we read this axiom, we get the impression  that it is a statement requiring proof, an impression reinforced by the fact that in the \emph{Elements},\index{Euclid!\emph{Elements}} the converse of this statement is a theorem (Proposition 27 of Book I)\footnote{Proposition 27 says: \emph{If a straight line falling on two straight lines make the
alternate angles equal to one another, the straight lines will be
parallel to one another.} \cite[p. 307]{Heath3}.} and by the fact that Euclid,\index{Euclid} having postponed the use of this axiom until the 29th proposition of Book I, tried to develop part of his geometry without it. Similar considerations had already been made and widely commented by Lambert's predecessors. Lambert recalls, in \S 9 of his memoir, that ``Proclus himself was troubled by Euclid's eleventh axiom, and required a proof of it, for the simple reason that \emph{its converse is provable}."
He mentions in \S 5 the conceptual difficulty (also raised many times before him) which arises from the fact that in using the parallel axiom, the lines must be extended to infinity. He considers that this use of infinity obfuscates the meaning of this axiom and makes it much less intuitive than the others. We may recall here that Lambert was also a philosopher and that, as such, the notion of infinity was not a mere word from the mathematical vocabulary but a concept that challenged him, as it challenged other philosophers of mathematics before him (like the neo-Platonic mathematician and philosopher Proclus, whom we have just quoted).

Lambert then discusses the ``representability" of the parallel axiom, declaring that if this axiom is not a theorem, then one should accept lines that are not \emph{straight} (Lambert emphasizes this word), but that can be asymptotes to each other. This cannot but  remind us of Lobachevsky's definition of parallelism:\index{parallel lines!definition!Lobachevsky}  two lines are parallel if they are asymptotes (in one direction), see \cite{L}.
Proclus, in his Commentaries on Book I of Euclid's \emph{Elements}, had\index{Euclid!\emph{Elements}} already made a similar statement; see e.g. his commentary on Definition 35  in \cite[p. 153 \&  seq.]{Proclus}.

Lambert in \S 8 quotes again Proclus,\index{Proclus of Lycia} saying that the latter did not advocate any change in the definition of parallelism but thought that one should be able to prove the parallel axiom without modifying any other axiom (\S 9). It appears from this passage that Lambert agrees with Proclus, who thought that Euclid had\index{Euclid} chosen to include the statement on parallels among his axioms because he had not found any proof.

  At the beginning of \S 11, Lambert brings up two possibilities in order
to solve the difficulties raised by the parallel axiom: either to deduce
it, as a theorem, from the other axioms, or to replace it by one or
several axioms that are simpler and clearer. He writes, at the end of
this introductive part: ``I have no doubt that Euclid had also thought
of including his eleventh axiom among the theorems". In the rest of
the memoir, Lambert will exploit the first path, of course without any success.

  \section{Part 2 of Lambert's memoir} \label{s:Lambert-Part-2}
Sections 12 to 26 constitute the second part of Lambert's memoir. In this part, Lambert presents a series of propositions that are valid in neutral geometry,\index{neutral geometry} that is, the geometry obtained from Euclidean\index{Euclidean geometry} geometry by keeping all the axioms and postulates, except that of parallels, and by adding no others.

Lambert begins by recalling (\S 12) that Euclid,\index{Euclid} until proposition 29 of the first book of the \emph{Elements},\index{Euclid!\emph{Elements}} did not use the parallel axiom, but that from this proposition onwards, he made extensive use of it.\footnote{In the commonly used versions of \emph{Elements}, the mention of parallel lines begins at proposition 27, not 29, of Book I, but as Lambert says, the use\index{parallel postulate} of the parallel postulate only begins at proposition 29. Propositions 27 and 28 state that if two straight lines make equal alternate internal (Proposition 27) or corresponding (proposition 28) angles with a third line that intersects them, then the lines are parallel.
These two propositions do not use the parallel postulate, but only the property  that the angle sum in a triangle is  $\leq 180^{\mathrm{o}}$, which is proved using Proposition 16 which does not use the parallel axiom. On the other hand, proposition 29 is the contrapositive of the parallel axiom (and is therefore equivalent to it).}
He notes that not only do the proofs of some propositions after the 29th use the parallel axiom, but that conversely, several of these propositions entail this axiom. He mentions as examples the proposition saying that the sum of the angles of a triangle is equal to two right angles (Proposition 32 of Book I) and another one saying that one straight line intersects another if and only if all the parallels to the first intersect the second, which is a proposition of Euclidean geometry\index{Euclidean geometry} that is not in the standard version of the \emph{Elements}.\index{Euclid!\emph{Elements}}

Lambert then announces that he will prove new propositions that do not use the parallel axiom. In fact, he proves some of the propositions he states, but for others, he does not succeed because, contrary to what he had in mind, they are not valid in neutral geometry:\index{neutral geometry} some are valid only under the negation of the parallel axiom, that is to say, in hyperbolic geometry, and others are valid only under the  the parallel axiom.

 \begin{figure}[htbp]
\centering
\includegraphics[width=1\linewidth]{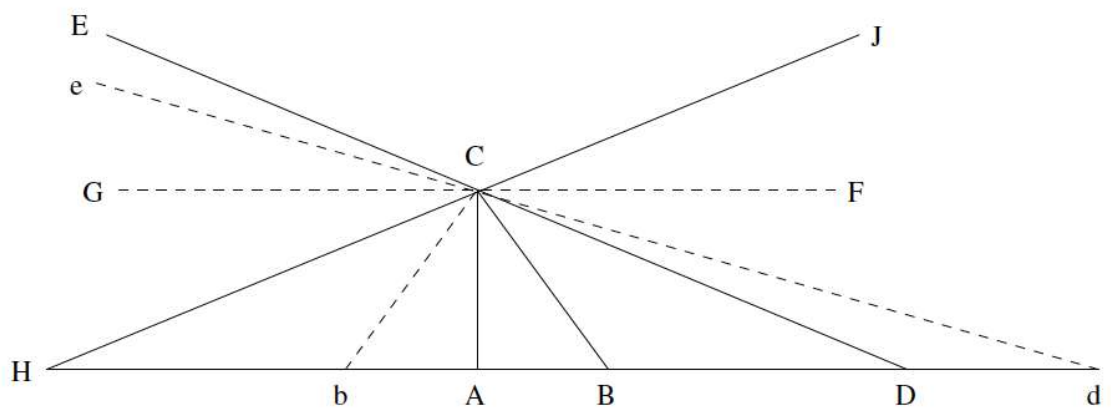}
\caption{\small  Figure for the Proposition in \S 13 of Lambert's memoir}
\label{fig:2-Lambert1-2}
\end{figure}

The first of these propositions is found in \S 13 of his memoir. He considers there a triangle $ABC$ which is right at $A$ and a line $ED$ that passes through $C$ and intersects the line $AB$ at $D$ (see Fig. \ref{fig:2-Lambert1-2}).
He shows that we have the following inequalities:
\[\widehat{ACD}< \widehat{ACB}+ \widehat{ABC} < \widehat{ACE}.\]

%
%
%
%
%

Lambert proves this proposition and then adds the following in \S 14: 
``This theorem now shows very precisely how far we are from Euclid's\index{Euclid} Proposition 17 concerning the determination of the sum of the three angles of a triangle"\footnote{Proposition 17 of Book I says: \emph{In any triangle two angles taken together in any manner
are less than two right angles.}} and he makes a long comment on his sentence. The interested reader is referred to \cite[p. 104 ff.]{2012-Lambert}

\begin{figure}[htbp]
\centering
\includegraphics[width=0.8\linewidth]{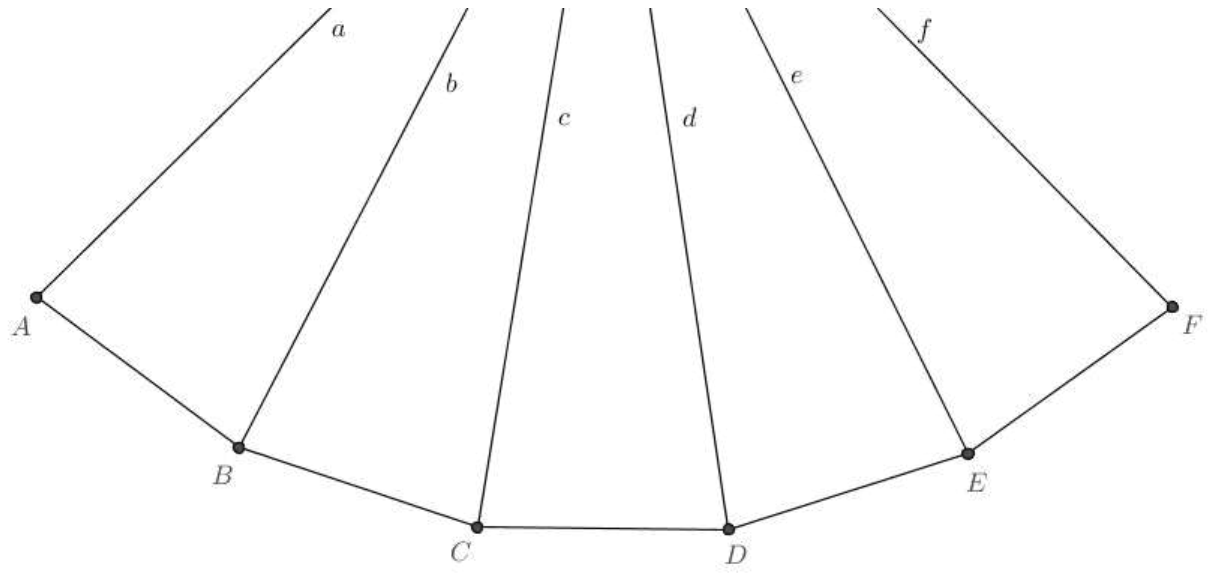}
\caption{\small Figure used for an argument in \S 15 of Lambert's memoir.}
\label{fig:2-Lambert2}
\end{figure}

In \S 15, Lambert says he wants to add another proposition, and he proves, by a simple geometric argument, that if the sum of the three angles in any triangle is constant, then this value is necessarily equal to two right angles. 
He then studies a property that can be stated in terms of regular polygons: Consider a sequence of points $A,B,C,D,E,F \ldots$ such that the consecutive segments $AB,BC,CD,DE,EF, \ldots$ are all equal and the angles at the points $B,C,D,E,F, \ldots$ are also all equal and strictly less than two right angles (Figure \ref{fig:2-Lambert2}). Consider the angle bisectors $Aa$, $Bb$, $Cc$, $Dd$, $Ee,\ldots$ of these angles. Then the sequence of points $A,B,C,$ $D,E,$ $F \ldots$ lies on a circle whose center is the common intersection of the bisectors $Aa$, $Bb$, $Cc$, $Dd$, $Ee,\ldots$

We know that this property does not hold in neutral geometry;\index{neutral geometry} in fact, it does not hold in hyperbolic geometry either, because the bisectors in question may not meet (consider consecutive equidistant points on a horocycle and join them by geodesic segments). This property is equivalent to the parallel axiom.

In \S 16, Lambert states the following property.

\begin{figure}[htbp]
\centering
\includegraphics[width=1\linewidth]{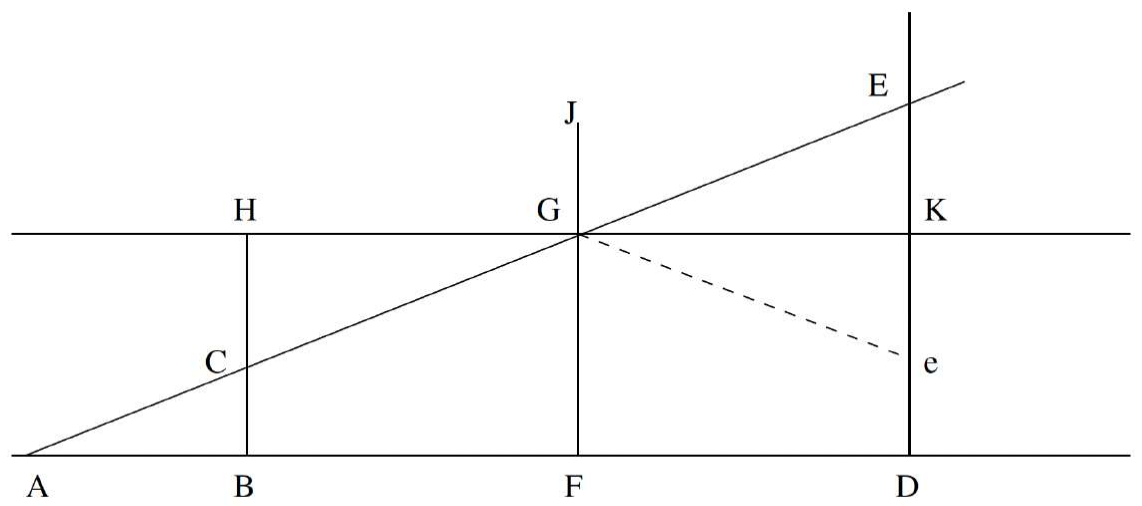}
\caption{\small Figure extracted from \S 16 of Lambert's memoir.}
\label{fig:2-Lambert24}
\end{figure}

\begin{quote}\small
Given a line $BD$ and a perpendicular $FG$ raised at a point $F$ on this line (see Figure \ref{fig:2-Lambert24} extracted from Lambert's memoir). Then from the point $G$ one can draw a line which meets $BD$ at a point $A$ such that the angle $\widehat{AGF}$ differs by as small an amount as one wishes from a right angle. 
\end{quote}

\begin{figure}[htbp]
\centering
\includegraphics[width=1\linewidth]{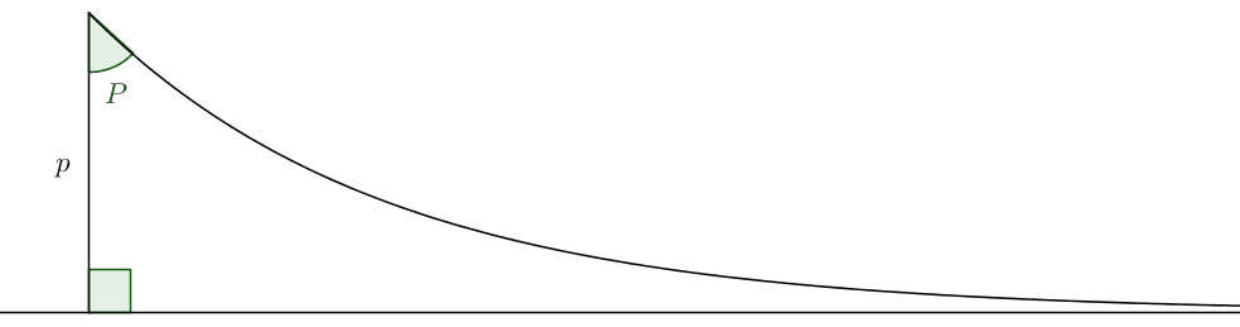}
\caption{\small $P$ is the angle of parallelism of the segment $p$}
\label{fig:2-Lambert4}
\end{figure}

After a preparation (\S 17, 18, 19), Lambert outlines three attempts to prove this property (\S 20, 21, 22) and then acknowledges that he does not succeed.
 We know that this property does not hold if the parallel axiom is not satisfied. Indeed, in neutral geometry,\index{neutral geometry} the angle $\widehat{AGF}$ has an upper limit, called the \emph{Lobachevsky parallelism angle} (Figure \ref{fig:2-Lambert4}),\index{Lobachevsky parallelism angle} whose value depends only on the length of the segment $GF$, and which is equal to a right angle if and only if the parallel axiom is satisfied. Under the negation of the parallel axiom, this upper limit is strictly smaller than a right angle. The reader may refer to the English edition of Lobachevsky's \emph{Pangeometry} \cite{L}, where Lobachevsky's angle of parallelism plays a paramount role.

\begin{figure}[htbp]
\centering
\includegraphics[width=1\linewidth]{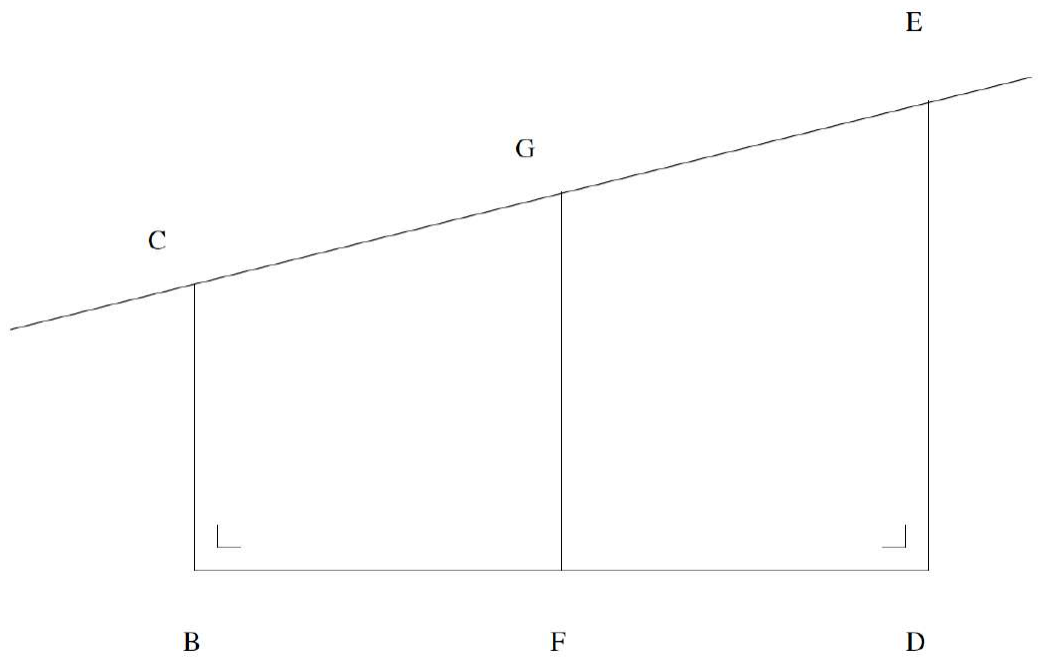}
\caption{\small Figure used in \S 23 of Lambert's memoir}
\label{fig:2-Lambert5}
\end{figure}

In \S 23, Lambert starts a study of quadrilaterals. He considers a quadrilateral $CBDE$ in which the angles $B$ and $D$ are right angles, and he calls $G$ the intersection point of $EC$ with the perpendicular to $BD$ raised at the point $F$, the midpoint of this segment (Figure \ref{fig:2-Lambert5}). He proves that we have $CB<DE$ if and only if the angle $\widehat{CGF}$ is acute.

In \S 23, 24 and 25, Lambert studies the monotonicity of the distance function between a point moving on a straight line and a second straight line. We know that this monotonicity requires the parallel axiom, and in fact, it equivalent to it.  Lambert eventually concludes (end of \S 25) that he cannot prove it without this axiom.

\begin{figure}[htbp]
\centering
\includegraphics[width=0.7\linewidth]{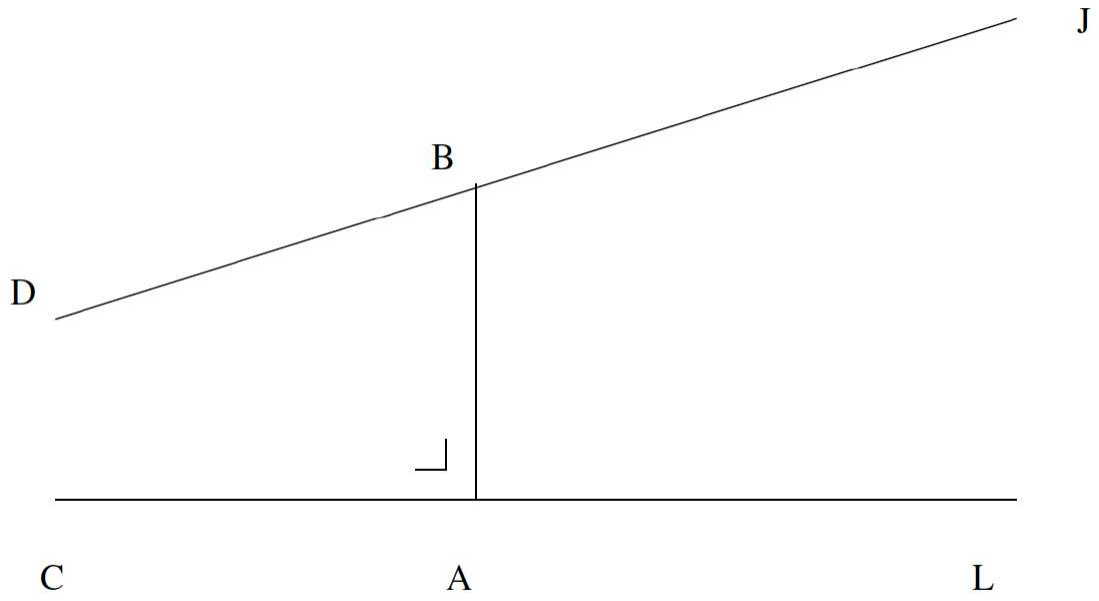}
\caption{\small Figure used in \S 26 of Lambert's memoir}
\label{fig:2-Lambert6}
\end{figure}

In \S 26, Lambert proves the following fact (without resorting to the parallel axiom).
Consider two lines $DJ$ and $CL$ and a point $B$ moving on $DJ$, from which we draw a perpendicular $BA$ onto $CL$ (Figure \ref{fig:2-Lambert6}). If the angle $\widehat{DBA}$ is always acute (regardless of the position of point $B$) then this angle is constant.

\section{Part 3 of Lambert's memoir} \label{s:Lambert-Part-3}

The third part of Lambert's memoir (\S 27 to 88) is the longest. It is titled \emph{Theory of parallel lines}, which is also the title of the whole memoir.
At the beginning of this part (\S 27), Lambert states that his intention is to present a theory that would find its place immediately after Proposition 28 of Book I of Euclid's \emph{Elements},\index{Euclid!\emph{Elements}} that is, that would be a continuation of the part of these \emph{Elements}  that does not need the parallel axiom. In this way, Proposition 29 and the ones that follow it, instead of resorting to the parallel axiom, would use Lambert's theory.

Lambert states that he considers that certain assertions made by his predecessors, such as the one on the equidistance of parallel lines, are not proved, and that he intends to prove them after having shown the difficulties contained in the statements of these propositions. He first summarizes the approach he will follow, whose central stage consists of successively examining three hypotheses, which will be that the fourth angle of a trirectangular quadrilateral is respectively right, obtuse or acute, and to prove that the second and third hypotheses contradict certain axioms of neutral geometry.\index{neutral geometry}
Trirectangular quadrilaterals\index{trirectangular quadrilateral}  were\index{quadrilateral!trirectangular} later given the name \emph{Lambert quadrilaterals}\index{quadrilateral!Lambert}\index{Lambert quadrilateral} or\index{Ibn al-Haytham--Lambert quadrilaterals} \emph{Ibn al-Haytham--Lambert quadrilateral}.\index{quadrilateral!Ibn al-Haytham--Lambert} Indeed, the Arab mathematician Ibn al-Haytham had already, more than seven centuries before Lambert, reduced his approach to the theory of parallels to the analysis of such quadrilaterals. In view of the importance of this topic, we have included, in  \S \ref{excursus:quadrilaterals} of this chapter, a detailed discussion of quadrilaterals in the theory of parallels.

\begin{figure}[htbp]
\centering
\includegraphics[width=0.8\linewidth]{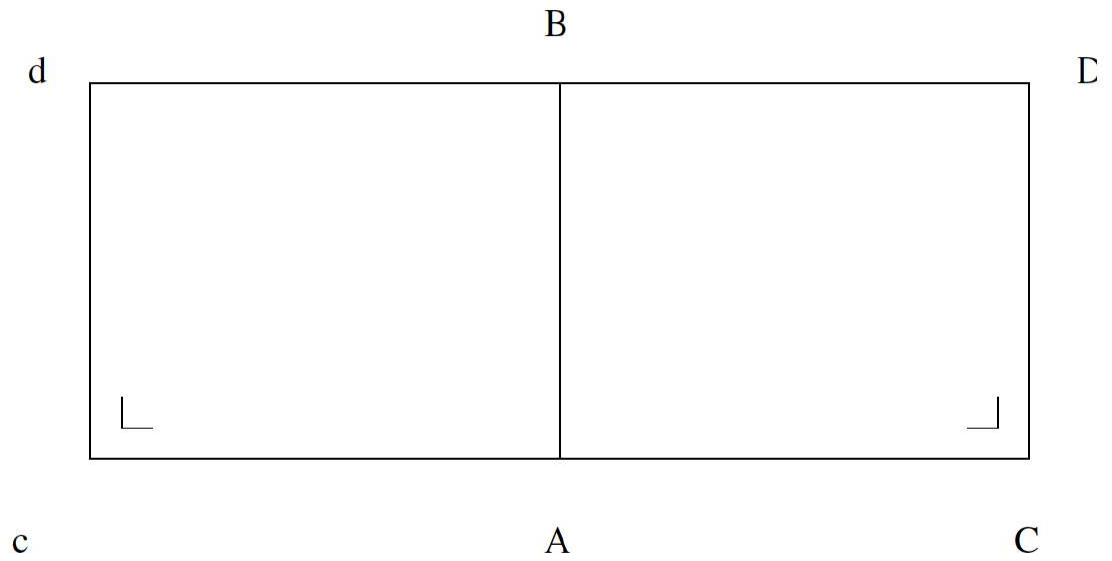}
\caption{\small The birectangular isosceles quadrilateral $cCDd$ is obtained by folding the trirectangular quadrilateral $ABCD$ along the side $AB$; the angles at $c$ and $C$ are right and we have $cd=CD$.}
\label{fig:2-Lambert7}
\end{figure}

In \S 28 to 38, Lambert recalls some properties of parallel lines that are contained in the \emph{Elements}\index{Euclid!\emph{Elements}} and some others that easily derive from them. He presents some simple properties of trirectangular quadrilaterals\index{trirectangular quadrilateral}  in neutral geometry.
In \S 39, he states again the three hypotheses relating to the trirectangular quadrilateral $ABCD$, represented in Figure \ref{fig:2-Lambert7}, in which the three angles $A$, $B$ and $C$ are right. This figure taken from Lambert's text also represents\index{quadrilateral!birectangular isosceles} the isosceles\index{quadrilateral!birectangular isosceles} birectangular quadrilateral\footnote{Birectangular isosceles quadrilaterals are also called \emph{Khayy\=am--Saccheri quadrilaterals}; see the review we make of quadrilaterals in \S \ref{excursus:quadrilaterals}.} $cCDd$\index{Khayy\=am--Saccheri quadrilateral} in which $cd=CD$, obtained by folding the quadrilateral $ABDC$ along the side $AB$.\footnote{Note that in this figure, as in the other figures, Lambert uses the following principle of symmetry in the choice of letters: $c$ and $d$ are the symmetrical points of $C$ and $D$ by the reflection with respect to the line $AB$. The use of a related notation is already found in the writings of Archimedes. The rule is practical and it indicates that if we exchange all the capital letters with the corresponding lowercase letters, the figure remains valid.}

 Lambert's three\index{trirectangular quadrilateral}  hypotheses are then the following:
\begin{enumerate}
\item $\widehat{BDC} = 90^{\mathrm{o}}$ (\emph{the right angle hypothesis})~;
\item $\widehat{BDC} > 90^{\mathrm{o}} $ (\emph{the obtuse angle hypothesis})~;
\item $\widehat{BDC} < 90^{\mathrm{o}} $ (\emph{the acute angle hypothesis}).
\end{enumerate}

The first hypothesis is examined in Sections 40 to 51, the second one in Sections 52 to 64, and the third one in Sections 65 to 88.

Under the first hypothesis, the quadrilateral $cCDd$ is a rectangle (i.e., a quadrilateral with four right angles; note that this implies that opposite sides are pairwise congruent). Lambert shows that in this case, for any\index{birectangular isosceles quadrilateral} isosceles\index{quadrilateral!birectangular isosceles} birectangular quadrilateral\index{Khayy\=am--Saccheri quadrilateral} (and not only the quadrilateral $cCDd$ on which we made our hypothesis), the perpendiculars drawn from any base to the opposite side are perpendicular to both sides and all have the same length (\S 41 and 42).

 Lambert then proves (\S 43 to 49) that any straight line $L$ passing through a point on one side of the quadrilateral and making an oblique angle with this side necessarily intersects the opposite side. His proof is based on a monotonicity property we state as follows: by placing side by side a sequence of rectangles congruent to $ABCD$ on a common base, the lengths of the various intersections of the straight line $L$ with the sides crossed in the successive quadrilaterals decrease by a quantity which is bounded from below (stated in \S 47). The result proved implies that the parallel axiom is valid under the hypothesis of the right angle. The aim of the rest of the memoir is to prove that hypotheses (2) and (3) are not compatible with Euclid's\index{Euclid} axioms deprived of the parallel axiom.

In Sections 48 to 51, Lambert tries to reverse the arguments  and he explains the difficulties encountered  in this process. In doing so, he proves a result, which he considers as remarkable, saying that if the first hypothesis holds in one rectangle then it holds in any rectangle. In fact, in doing this, Lambert proves a series of results that the 17th-18th century mathematician Giovanni Girolamo\index{Saccheri, Giovanni Girolamo}  Saccheri\footnote{\label{n:Saccheri} Giovanni Gerolamo Saccheri (1667-1733) was a grammar teacher before becoming interested in mathematics, probably under the influence of his teacher and friend the Jesuit priest Tommaso Ceva\index{Ceva, Tommaso}, brother of the famous mathematician Giovanni Ceva\index{Ceva, Giovanni}. Saccheri entered the Jesuit order in 1685 and was ordained in 1694. He was appointed professor of philosophy at the University of Turin that same year, and three years later at the University of Pavia, where he taught philosophy, theology and mathematics. He took a close interest in Euclid's \emph{Elements}, following the advice of his first teacher, Ceva. He held this position until his death. Like other great mathematicians, he tackled the question of parallels, and he attempted --- a victim, like all his predecessors, of the prejudices of his time --- to prove that the parallel postulate is a theorem.} had obtained in his treatise \cite{Saccheri}, where he makes three similar hypotheses about triangles. Saccheri proved that if in a particular triangle one of the three hypotheses is valid, then it is also valid in any triangle. The result was later called \emph{Saccheri's theorem} by Bonola in \cite{Bonola-2}.
 
\begin{figure}[htbp]
\centering
\includegraphics[width=1\linewidth]{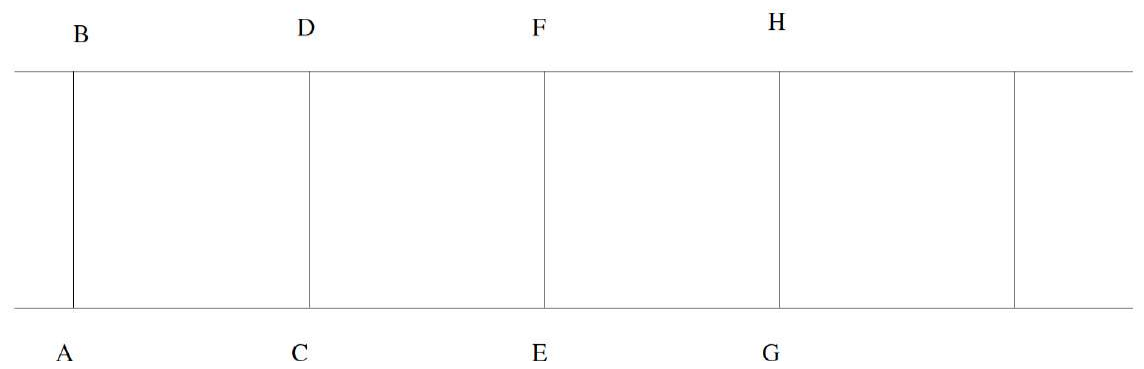}
\caption{\small{At the points $A$ and $B$ and at points $C,E,G,...$ the angles are right angles. Under the second hypothesis, the segments $CD, EF, GH,...$ decrease more rapidly than linearly with respect to the distances from points $C,E,G,...$ to point $A$. Similarly, the angles at $D,F,H,...$ become increasingly obtuse.}}
\label{fig:2-Lambert8}
\end{figure}

 In Sections 52 to 54, Lambert proves some elementary properties of the trirectangular quadrilateral\index{trirectangular quadrilateral} under\index{Ibn al-Haytham--Lambert quadrilateral} the obtuse angle hypothesis.  He shows in particular that each of the two sides surrounded by right angles is longer than the opposite side. He then considers (Sections 55 to 57) perpendiculars drawn onto the line containing a side adjacent to two right angles from a sequence of points located on the line containing the opposite side and which are outside the rectangle, on the side of the acute angle and extending further than this acute angle (Figure \ref{fig:2-Lambert8}).

He proves that the sequence of lengths of these perpendicular segments
is monotonically decreasing, and that the difference between one segment and the previous one is smaller the further apart the points are. It may be recalled here that such monotonicity properties for lengths in spherical geometry are well known, and can be found in the \emph{Spherics}  of Theodosius of Tripoli (c. 60-90 BC) and in those of Menelaus of Alexandria (c. 70-140 AD).\footnote{The \emph{Spherics} of Theodosius and those of Menelaus are two comprehensive treatises on spherical geometry, that constitute the most important works on this subject that have come down to us from Ancient Greece. It should be noted that both treatises (and especially the second) are very rich from the mathematical point of view. Greek manuscripts of Theodosius' treatise survive, and were the bases of modern editions, but all 
the Greek versions of Menelaus treatise are lost, and this treatise has come down to us only in the form of Arabic translations. The most faithful edition of Theodosius' \emph{Spherics} remains the one in French, by P. Ver Eecke \cite{Ver-Eecke-Theodosius}, and it contains valuable notes. Concerning  Menelaus' \emph{Spherics}, there is a German translation by M. Krause, based on the Arabic text of Ibn `Ir\=aq \cite{Krause}, and a recent critical edition with an English translation based on the Arabic translations of al-M\=ah\=an\=\i \ and al-Haraw\=\i,   by R. Rashed and A. Papadopoulos \cite{RR}. The treatises of Theodosius and Menelaus are very different in their approach and in the techniques used. The proofs of the propositions of Theodosius \emph{Spherics} are based on the Euclidean geometry of the ambient space of the sphere  considered as embedded in three-dimensional space, whereas those of Menelaus, except for a very small number of propositions, are intrinsic. The monotonicity theorems on lengths to which we refer here are Propositions 6 to 10 of Book III of Theodosius' \emph{Spherics}, and there are a large number of them in Menelaus' \emph{Spherics}, including Propositions 45-61 and 81-91 of al-M\=ah\=an\=\i /al-Haraw\=\i's version \cite{RR}.}

After this non-linear monotonicity property for segment lengths, Lambert proves (\S 58 to 60) a monotonicity property for the angles that the perpendiculars make with the line on which the point from which they are issued is located. These angles become more and more obtuse (Figure \ref{fig:2-Lambert8} again). Monotonicity properties of angles are also known in spherical geometry.\footnote{See Propositions 84-86  of Menelaus' \emph{Spherics} \cite{RR}.}

Lambert deduces from this that the lengths of the perpendicular segments cannot become asymptotically shorter  (\S 62) but that these lengths become zero at a certain point, and consequently the line containing the points from which the perpendiculars are drawn must intersect the line containing the opposite side of the quadrilateral. This contradicts the Euclidean\index{Euclid} axiom\index{Euclidean geometry} which says that two lines cannot intersect at more than one point (\S 64). The second hypothesis is therefore invalid.

 Lambert was aware of the fact that his second hypothesis is valid on the sphere, but he would only state this in \S 82, when analyzing his third hypothesis, and when he would state (in comparison) that the latter is realized on a sphere of imaginary radius.

Starting in \S 65, Lambert examines the third hypothesis, that of the acute angle. He begins with considering a\index{trirectangular quadrilateral}  trirectangular\index{quadrilateral!trirectangular} quadrilateral (whose fourth angle is therefore acute) and proves (\S 66 and 67) that each of the two sides adjacent to the acute angle is greater than the side opposite to it. He then shows (\S 68) that each perpendicular drawn from the extension of one side of the acute angle to the opposite side is shorter than a perpendicular drawn from a point on the same line but further away. Just as under the hypothesis of the obtuse angle, Lambert studies (\S 69) the monotonicity of the angles formed by such lines. These angles become more and more acute.

Then, and as in the previous case, Lambert shows (\S 70) that the lengths of the perpendiculars not only become larger and larger, but that they can become larger than any value given in advance, whereas under the obtuse angle hypothesis, such strict monotonicity was sufficient to contradict one of the axioms of the \emph{Elements}.\index{Euclid!\emph{Elements}} From this point on, Lambert rediscovers a property that Saccheri\index{Saccheri, Giovanni Girolamo} had already emphasized and which is specific to hyperbolic geometry (which, let us recall, was for Lambert, as for Saccheri, only hypothetical): the existence of disjoint lines having a common perpendicular and diverging from each other on each side of this perpendicular. Lambert is struck by this property which makes him imagine lines that do not resemble the usual lines. However, he refuses to be convinced by these representations.

Note that if Lambert had considered, like many of his predecessors, that parallelism implies equidistance between lines, or that equidistance defines parallelism, the work he has undertaken here would already be complete. But Lambert seeks another argument.

\begin{figure}[htbp]
\centering
\includegraphics[width=1\linewidth]{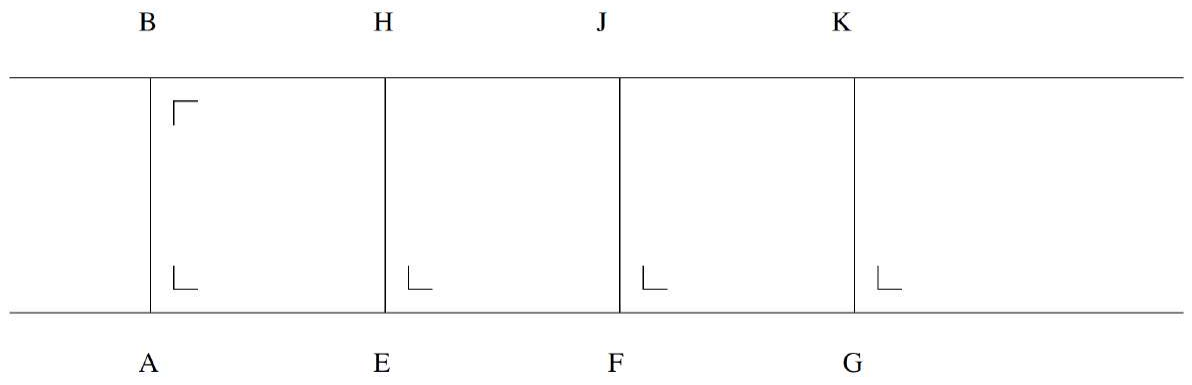}
\caption{\small Figure used in \S 72 and 83 of Lambert's memoir.}
\label{fig:2-Lambert9}
\end{figure}

In \S 72, Lambert discusses two properties that could hold under the acute angle hypothesis. They are related to Figure \ref{fig:2-Lambert9}, where we have two lines AE and BH with a common perpendicular AB and where the angles at the sequence\index{quadrilateral!trirectangular} of points $E, F, G, \ldots$ are right angles. The quadrilaterals $BAEH, BAFJ, BAGK, \ldots$ are trirectangular.\index{trirectangular quadrilateral}  The first possibility  mentioned by Lambert here is that it is conceivable that when the points $E,F,G,\ldots$ are sufficiently far from the point $A$, the perpendicular raised at these points does not meet the line $BH$. The second possibility is that the angle at $H, J, K,\ldots$ which, as we have already seen, decreases when the point moves away from $A$, becomes eventually zero, not only at the limit, but for a given point on $BH$.

 In \S 73 and \S 74, Lambert proves that under the acute angle hypothesis, the sum of the three angles of each triangle is $< 180^{\mathrm{o}}$.

In \S 76 and \S 77, he shows that in an equilateral triangle $ABC$, if $F$ is the midpoint of $BC$ and $D$ the intersection point of the medians of this triangle, then we have $DF < \frac{1}{3}  AF$. It is easy to see that in neutral geometry,\index{neutral geometry} the three medians of an equilateral triangle are concurrent. The case of an arbitrary triangle is more difficult to handle, but Lambert does not discuss it. He notes at the same time that under the obtuse angle hypothesis, we would have $DF> \frac{1}{3} AF$ (and it is well known that $DF= \frac{1}{3} AF$ under the right angle hypothesis).

Lambert notes in \S 79 that there is a series of results, analogous to the one just stated and to several others that we saw earlier (monotonicity of lengths, angles, etc.), valid under the hypothesis of the acute angle and obtained by reversing the direction of the conclusions of the corresponding results obtained under the hypothesis of the obtuse angle. He notes, however, that, despite this resemblance in the results, the third hypothesis (that of the acute angle) does not lead as easily as the second one (that of the obtuse angle) to a contradiction.

According to Lambert, one of the most important consequences of the third hypothesis is the existence of an absolute measure for lengths, areas, and volumes. He proceeds (\S 80) as follows.
We start from the fact that there is an absolute measure for\index{Gauss, Carl Friedrich} angles.\footnote{Declaring that the total angle at a point is equal to four right angles does not depend on any pre-chosen standard measure. It is in this sense that we have an absolute measure for angles. Using the fact that the lengths of the sides of a triangle are completely determined by their angles, we obtain an absolute measure for lengths. Several years after Lambert, Gauss  had also noted this fact. In a letter to Gerling, dated April 11, 1816 (see \cite[Vol. VIII, p. 168]{Gauss}), Gauss wrote: ``It would have been desirable if Euclidean geometry were not true, because we would have a universal measure a priori. We could use the side of an equilateral triangle with angles 59${}^{\mathrm{o}}$59'59'',9999.....as a unit of length." [z.B. k\"onnte man als Raumeinheit die Seite desjenigen gleichseitigen Dreiecks annehmen, dessen Winkel = 59${}^{\mathrm{o}}$59'59'',9999..... ]. We recall that in hyperbolic geometry, the common value of the angles in an equilateral triangle is always strictly less than 60${}^{\mathrm{o}}$, and this is why Gauss chose such a value.} From this, Lambert notes that there is an absolute measure for the areas of isosceles trirectangular quadrilaterals,\index{trirectangular quadrilateral}  i.e., those\index{quadrilateral!trirectangular} for which the two sides not adjacent to the acute angle have the same length. Such a quadrilateral is completely determined by the value of the acute angle, and the area of this quadrilateral depends solely on the value of this angle. Thus, the value of this angle is an absolute measure of the area of the quadrilateral. Lambert states that ``this consequence has such a force of attraction that it would easily force us to wish that the third hypothesis were true!", but he notes that if this were true, then ``trigonometric tables would become interminably long,\footnote{In non-Euclidean geometry, in a right triangle, the ratio of the length of one side to that of the hypotenuse does not only depend on the angle between these two sides; it also depends on the size of the triangle. It is in this sense that we understand the statement that the trigonometry tables, in hyperbolic geometry, are of infinite length.} the similarity and proportionality of the figures would have completely disappeared~;
no figure could be represented other than in its absolute size."

In \S 81 and 82, Lambert studies the area function of triangles. He shows that the area of a triangle can be defined as its \emph{angular defect}, that is, the difference between $180^{\mathrm{o}} $ and the sum of the angles. He also points out that under the second hypothesis (that of the obtuse angle), the area of a triangle can be defined as its \emph{angular excess},\index{Cavalieri, Bonaventura} that is, the excess over $180^{\mathrm{o}} $ of the sum of the three angles of the triangle\footnote{This theorem is attributed to Albert Girard (1595-1632), \index{Girard, Albert} who stated it in his \emph{New invention in algebra} (1629). Lagrange, \index{Lagrange@de Lagrange, Joseph-Louis} mentions, in his memoir \emph{Solution of some problems relating to spherical triangles with a complete analysis of these triangles} (1800) translated in Chapter 18 of the present volume, that Girard's proof is not sufficient and he considers that this theorem should rather be attributed to Cavalieri, who gave a complete proof in his \emph{Directorium generale uranometricum} (General guide to celestial measurements), Bologna, 1632.}. Lambert states that ``we should almost conclude that the third hypothesis occurs for an imaginary sphere."

This remark about the imaginary sphere is interesting in several respects, and can be considered prophetic in a certain sense. Indeed, the area of a spherical triangle, on a sphere of radius $r$, is equal to the angular excess of this triangle, that is, to the quantity $r^2(\alpha+\beta+\gamma-\pi)$ where $\alpha, \beta, \gamma$ are the values of the three angles of the triangle. If instead of the radius $r$ we take an \emph{imaginary radius} and multiply this value by the imaginary number $\sqrt{-1}$, we obtain, for the area of a triangle, $r^2 (\alpha+\beta+\gamma-\pi) =  (\sqrt{-1})^2 (\alpha+\beta+\gamma-\pi) = -(\alpha +\beta +\gamma -\pi) = \pi-\alpha-\beta-\gamma$, which is the area of a triangle in the hyperbolic plane whose angle values are $\alpha, \beta, \gamma$.
It may be noted here that Lambert, in his \emph{Trigonometric Observations} \cite{Lambert-Observations}, highlights the analogy between the trigonometric functions sine and cosine and the hyperbolic sine and cosine functions, noting that the hyperbolic cosine of an arc can be considered as the cosine of an imaginary arc, thus developing a kind of trigonometry in which the sides of triangles are considered as imaginary. But he does not seem to have made the connection between this and a geometry where the sum of the angles in a triangle is smaller than $180^{\mathrm{o}} $, even though he was very close to these ideas.

It was not until Franz Adolph Taurinus\index{Taurinus, Franz Adolph} (1794-1874) and then Lobachevsky\index{Lobachevsky, Nikolai Ivanovich} that the statement on the imaginary sphere in relation with the area of non-Euclidean\index{Euclid} triangles was formulated more explicitly. Taurinus is the author of another {\it Theorie der Parallellinien} also reproduced in the volume by St\"ackel and Engel which contains Lambert's text on the theory of parallels, in which Taurinus also attempts to prove Euclid's parallel axiom from the other axioms. Also seeking a contradiction in the negation of this axiom, Taurinus ends up, like Lambert, developing the elements of a non-Euclidean geometry. We owe to Taurinus one of the first attempts to establish a parallel between the trigonometric formulae of the three geometries: spherical, Euclidean\index{Euclidean geometry} and hyperbolic (the existence of the latter having been, for him, as for Lambert, purely hypothetical).
 In his memoir, Taurinus obtained fundamental formulae of hyperbolic trigonometry by working on a sphere of imaginary radius, calling the geometry obtained in this way a {\it logarithmic-spherical geometry} (\emph{Logarithmisch-sph\"arischen Geometrie}). He noted a formal transition between the formulae of hyperbolic geometry and those of spherical geometry which consists in replacing in the spherical formulae certain trigonometric functions ($\sin$ and $\cos$) by hyperbolic functions ($\sinh$ and $\cosh$). It can also be noted that, in the same way, Lobachevsky,\index{Lobachevsky, Nikolai Ivanovich} in his {\it Elements of geometry} (1829-1830) and in his {\it Geometrical research on the theory of parallels} (1840), highlights a passage between the formulae of spherical trigonometry and those of hyperbolic trigonometry which consists of replacing, in the formulae of spherical trigonometry, the lengths of the sides $a,b,c$ of a triangle by the imaginary quantities $a\sqrt{-1},b\sqrt{-1},c\sqrt{-1}$. Note also that the transition from spherical trigonometry formulae to hyperbolic trigonometry formulae involved the replacement of certain trigonometric functions with hyperbolic functions, and that Lambert was one of the first to systematically use hyperbolic functions. Finally, note that Lambert was among the mathematicians who developed the geometric foundations of the theory of complex numbers.

Let us now return to Lambert's \emph{Theory of the parallel lines}.

In \S 83, Lambert proves a property he had already mentioned
in \S 72, which states (using the notation in Figure \ref{fig:2-Lambert9}) that the perpendiculars raised at points $E,F,G,\ldots$ to the line $AE$ do not necessarily meet the line $BH$, when the point where the perpendicular is raised is far enough from $A$.

The last paragraphs of the memoir (\S 84-88) contain sketches of ideas for proving that the acute angle hypothesis leads to a contradiction, but Lambert acknowledges that most of them are unsuccessful. In the last paragraph, he mentions a new attempt, in passing over several details, and the text does not allow us to know whether he himself was convinced of this or not.

\section{Excursus: On the theory of parallel lines in history}\label{excursus:parallel-lines}

\subsection{Greek Antiquity}
We recalled in \S \ref{s:Lambert-Part-1} that Lambert refers to Proclus regarding the fifth postulate.
Let us quote\index{Proclus of Lycia} Proclus' passage in question (Proclus 192, \cite[p. 150-151]{Proclus-Morrow}):

 \begin{quote}\small
 
[The fifth postulate] ought to be struck from the postulates altogether. For it is a theorem --- one that invites many questions, which Ptolemy proposed to resolve in one of his books --- and requires for its proofs a number of definitions as well as theorems. And its converse is proved by Euclid himself as a theorem. But perhaps some persons might mistakenly think that this proposition deserves to be ranked among the postulates on the ground that the angles being less than two right angles makes us at once believe in the convergence and intersection of the straight lines. To them Geminus\index{Geminus of Rhodes} has given the proper answer when he said that we have learned from the very founders of this science not to pay attention to plausible imaginings in determining what propositions are to
be accepted in geometry. Aristotle likewise says that to accept probable reasoning\footnote{The expression ``probable reasoning" here refers to the fact that lines would ``eventually" intersect, if they are sufficiently extended.} from a geometer is like demanding proofs from a rhetorician.\footnote{Aristotle, \emph{Nicomachean Ethics}, 1094b26f.} And Simmias is made by Plato to say, ``I am aware that those who make proofs out of probabilities are impostors."\footnote{Plato, \emph{Phaedo} 92d.} So here, although the statement that the straight lines converge when the right angles decrease is true and necessary, yet the conclusion that because they converge more as they are extended farther they will meet at some time is plausible, but not necessary, in the absence of an argument proving that this is true of straight lines. That there are lines that approach each other indefinitely but never meet seems implausible and paradoxical, yet it is nevertheless true and has been ascertained for other species of lines. May not this, then, be possible for straight lines as for those other lines? Until we have firmly proved that they meet, what is said about other lines strips our imagination of its plausibility. And although the arguments against the intersection of these lines may contain much that surprises us, should we not all the more refuse to admit into our tradition this unreasoned appeal to probability?
 \end{quote}

Let us say a few more words on Proclus' own\index{Proclus of Lycia} contribution. We have repeatedly pointed out that his work contains valuable information on the history of Greek geometry. It is not certain that Proclus obtained original results, but a ``proof"
of the parallel axiom in his \emph{Commentaries on Book I of Euclid} \index{Euclid} is attributed to him. Not being satisfied with the proofs of this axiom given by his predecessors, Proclus gave his own proof, based on the hypothesis that if two coplanar lines are disjoint, their distance is bounded, a hypothesis which, as already noted, is equivalent to the parallel axiom. Based on work by Geminus,\index{Geminus of Rhodes} Proclus conceived the possibility of a geometry involving lines that are not necessarily straight but which can be asymptotes. We have already quoted a passage in which Proclus discusses this question; see also his commentary on Definition 35, p. 153 et seq. of \cite{Proclus}. All this has led some commentators to consider Proclus as one of the precursors of hyperbolic geometry.\index{hyperbolic geometry}

We also noted that doubts concerning the parallel axiom had already been raised before Euclid. Aristotle\index{Aristotle} (384-322 BC) discusses the theory of parallels and the difficulties posed by this theory in the {\it First Analytics} (II, 16, 65a 4;II, 17, 66a 11-15) and the {\it Second Analytics} (II, 15, 65a 4-7).  
He declares there that there is a vicious circle in the theory of parallels, that is, something that is impossible to prove, ``because of a difficulty in the definition of parallels".\index{parallel lines!definition!Euclid}   Archimedes (ca. 287-212 BC) wrote a treatise, \emph{On parallel lines}, which does not survive but which is quoted by the Arab bibliographer al-Nad\=\i m in his \emph{Book of Bibliography} (Kit\=ab Fihrist).
 
We know from Proclus that Heron of Alexandria\index{Heron of Alexandria} (10-70 AD) admitted only three axioms; they are quoted in \cite[p. 173]{Proclus}. Proclus writes: ``However, we should not reduce the number of axioms to a minimum like Heron, who only set three". In the same passage, Proclus\index{Proclus of Lycia} tells us that\index{Pappus of Alexandria} Pappus of Alexandria (ca. 290-350) added other axioms to those of Heron.

Ptolemy\index{Ptolemy of Alexandria} (90-168) wrote a treatise on parallels, titled {\it On the intersection of straight lines}, in which he gave a proof of the parallel axiom, based on the assumption that two parallel lines intersecting a third make equal alternate interior angles with it. This assertion is equivalent to the parallel axiom. He attempted to prove Euclid's proposition I.29 without using the parallel axiom. This work has not reached us, but this proof is reported by Proclus.\footnote{See \cite{Proclus} p. 312 et seq.}

  Attempts  by the Neoplatonic philosopher Simplicius\index{Simplicius} (c. 490 - c. 560) to prove the parallel postulate\index{parallel postulate} have come down to us in the form of extensive excerpts in the writings of the 10th-century\index{Nayrizi@Al-Nayr\=\i z\=\i} mathematician Al-Nayr\=\i z\=\i . This brings us to the Arab mathematicians.

\subsection{The Arabs}

In the eighth century AD, an Arabic-language science was born. It quickly reached a high degree of maturity and remained alive for about six centuries, with some very talented mathematicians, several of whom involved in the problem of parallels. We shall review some of their works on this subject. To tell the story in a few simple words, the Byzantines,
who had ordered the closure of the Philosophical schools of the Empire, had at the same time condemned Hellenistic mathematics, which was taught in the Philosophical schools. Taking advantage of this opportunity, from the eighth century onwards, several Arab missions were sent to the Byzantine court with the aim of acquiring copies of ancient mathematical manuscripts. Because of their significant interest in the Greek's mathematics and philosophy, the Arabs became their natural heirs. They studied Aristotle,\index{Aristotle} for whom they had great respect, and whom they often called ``the philosopher" or ``the wise".  They extensively translated and commented his works along with those of Greek mathematicians, and the problems raised by the question of parallels naturally emerged in these works.

`Abb\=as Ibn Sa`\=\i d al-Jawhar\=\i\   (ca. 800-860),\index{Jawhari@al-Jawhar\=\i } who worked in Bagdad and Damascus, is one of the first known Arab mathematicians who studied the parallel problem.\index{Tusi@al-\d{T}\=us\=\i , Na\d{s}\=\i r al-D\=\i n}  Na\d{s}\=\i r al-D\=\i n al-\d{T}\=us\=\i , in his \emph{Rectification of the book of the Elements}, describes an attempt of a proof of the fifth postulate by al-Jawhar\=\i ,\index{al-Jawhar\=\i} based on the fact that through any point inside an angle whose value is strictly between $0$ and $180^{\mathrm{o}}$, we can draw a line that intersects both sides. We know that this property is equivalent to the parallel axiom and therefore does not follow from the other axioms. Al-\d{T}\=us\=\i \ described new propositions added by Al-Jawhar\=\i\ to Euclid's \emph{Elements}, see Rosenfeld \cite{Rosenfeld} p. 47.

We also mention Th\=abit Ibn Qurra\index{Ibn Qurra, Th\=abit} (826-901), an Arab mathematician whose mother tongue was Syriac and who worked in Mesopotamia (Harr\=an, then Baghdad). He translated treatises of Apollonius\index{Apollonius of Perga}, including several books of the \emph{Conics}, as well as works of Archimedes, Euclid and Ptolemy.\index{Euclid} He wrote at least two original treatises on parallels, one of them titled {\it On the Proof of Euclid's Famous Postulate} and another one titled {\it Book on the fact that if a straight line falls on two straight lines and makes with them and on the same side less than two right angles, the two straight lines, if extended on that side, meet}. Each of these treatises contains an (attempted) proof of the parallel axiom. The first proof is based on the assumption that if two lines diverge on one side, they converge on the other, and the second is based on the existence of equidistant lines. We know that both hypotheses are equivalent to the parallel axiom. The two treatises of  Th\=abit Ibn Qurra have come down to us. The reader is referred to the paper by Rashed and Houzel \cite{RH2005} which contains a French translation of these papers as well as a commentary and references to previous works.  The reader may also consult Rosenfeld \cite{Rosenfeld}. According to Rashed and Houzel \cite{RH2005}, Ibn Qurra was the first to give a rigorous proof of the symmetry of the equidistance relation (i.e., if all points of a line $d$ are equidistant from a line $d'$, then all points of $d'$ are equidistant from $d$), once the existence of equidistant lines is admitted (which Th\=abit Ibn Qurra\index{Ibn Qurra, Th\=abit}  attempted to prove in the second treatise we just mentioned).

We also mentioned  Al-Nayr\=\i z\=\i,\index{Nayrizi@Al-Nayr\=\i z\=\i}  the 10th-century mathematician who wrote a \emph{Treatise on the proof of Euclid's famous postulate}. He used the notion of equidistance as a definition of parallelism. Al-Nayr\=\i z\=\i \ quotes the Arabic translation of a text by Simplicius in which the latter reports on a proof\index{Agh\=anis} of Agh\=anis,\footnote{Agh\=anis is probably an Arabized name of a fifth-century Greek mathematician, who knew Simplicius.} having also ``proved" the parallel axiom using the notion of equidistance and having further assumed that this notion is symmetric.

 Other mathematicians of the 10th century include Yu\d{h}an\-n\=a Ibn Y\=usuf al-Quss\index{Ibn Y\=usuf al-Quss} (died c. 980), an Arab Christian priest who wrote a book dedicated to the Emir of Aleppo Sayf al-Dawla titled \emph{On the intersection of two straight lines rising from a straight line and making angles with it less than two right angles}, and Ab\=u `Al\=\i \ Ibn S\=\i n\=a \index{Ibn S\=\i n\=a} (981-1037), the famous philosopher, physicist, mathematician, physician, astronomer and music theorist, known in the Western literature under the name Avicenna. Ibn S\=\i n\=a's famous treatise, \emph{The book of science}, contains a chapter on geometry, in which parallel lines are defined using the concept of equidistance.

\begin{figure}[htbp]
\centering
\includegraphics[width=1\linewidth]{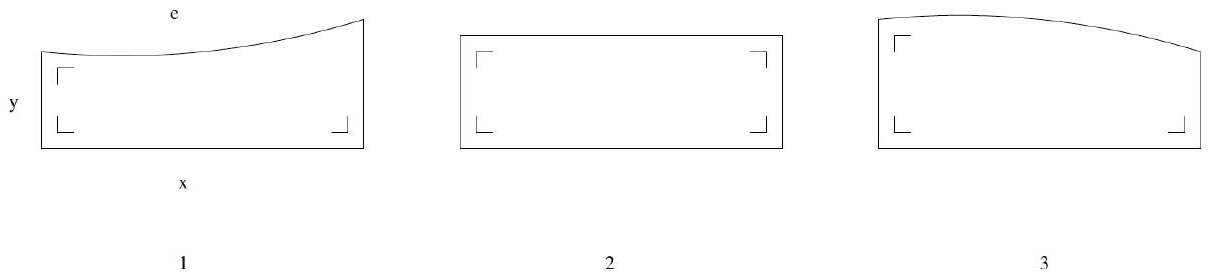}
\caption{\small Trirectangular quadrilaterals. In the figure on the left (hyperbolic geometry), the fourth angle is acute; in the middle figure it is right (Euclidean geometry), and in the right figure it is obtuse (spherical geometry).}
\label{fig:1trois3}
\end{figure}

Then one has to mention the great mathematician Ab\=u  `Al\=\i \ Ibn al-Haytham\index{Ibn al-Haytham} (d. after 1040), known in Latin as Alhazen.\footnote{Born in Basra and working in Cairo, Ibn al-Haytham was an engineer, physicist, and astronomer. He is sometimes called the father of modern optics. His treatise on optics has long been known in the West thanks to Gerard of Cremona's translation. A crater on the moon bears his name.} He\index{Euclid!\emph{Elements}!Arabic manuscripts} is the author of several treatises on geometry, including the \emph{Book explaining the postulates of Euclid's book on the Elements}\index{Euclid!\emph{Elements}} and the {\it Book on the resolution of doubts raised by Book I of Euclid's Elements}. In the first treatise, Ibn al-Haytham calls ``parallel to a given line" the locus of the end of a segment that moves perpendicularly to the given line while maintaining a constant length.
He proves\index{trirectangular quadrilateral} that\index{quadrilateral!trirectangular} this locus\index{Ibn al-Haytham--Lambert quadrilateral} is a straight line. His ``proof" of the parallel axiom is based on the consideration of a quadrilateral with three right angles,\footnote{Pont and Rosenfeld call such a quadrilateral, which Ibn Qurra had already considered, \emph{Ibn al-Haytham--Lambert quadrilateral}, see \cite[p. 171]{Pont} and \cite[p. 65]{Rosenfeld}. We use both terminologies.} see Figure \ref{fig:1trois3}.  

Ibn al-Haytham\index{Ibn al-Haytham} 
 analyzes the three cases.\index{trirectangular quadrilateral} Like\index{quadrilateral!trirectangular} Lambert in his treatise, he ``proves" that only the first hypothesis does not lead to a contradiction. This analysis is reviewed by Youschkevitch in \cite [p. 117]
{Youschkevitch}, Rosenfeld \cite[p. 59]{Rosenfeld} and Pont \cite[p. 169]{Pont}. 
With this work, Ibn al-Haytham was a direct precursor of Lambert, who followed exactly the same approach (except for the definition of parallels). We refer to \cite[vol. V]{Rashed1} for the translation and commentaries of the mathematical works of Ibn al-Haytham.

  Next, we must mention `Umar al-Khayy\=am\index{Khayy\=am@al-Khayy\=am, `Umar}, the mathematician born in Nay\-sh\=ab\=ur in 1048 and who died in the same city in 1131. He was also a philosopher and was interested in the philosophical
  foundations of mathematics. He studied the works of Ibn al-Haytham. According to Rashed (Introduction to \cite{Rashed0}), Khayy\=am developed a geometric theory of equations that can be considered as the beginning of algebraic geometry. Khayy\=am\index{Euclid} completed his \emph{Commentary on the difficulties of certain postulates in Euclid's work} in 1077. We have already quoted a passage from this work in which the author points out that his predecessors did not correctly address the question of parallels, claiming that they have neglected certain principles that Aristotle had stated. 
  
  Khayyam's work consists of three parts. The first, which interests us here, concerns the theory of parallels. The second concerns the theory of proportions, and the third concerns compounded ratios.
 
In the first part of this book, Khayy\=am studies (Proposition 3) quadrilaterals with two right angles adjacent to one side, the other two angles being equal. These quadrilaterals are the so-called \emph{birectangular isosceles rectangles} or \emph{Khayy\=am--Saccheri quadrilaterals},\index{Khayy\=am--Saccheri quadrilateral} and\index{birectangular isosceles quadrilateral} we\index{quadrilateral!birectangular isosceles} shall discuss them in more detail in \S \ref{excursus:quadrilaterals}. The three quadrilaterals in Figure \ref{fig:1trois3} represent the three cases of trirectangular quadrilaterals.\index{trirectangular quadrilateral} In the\index{quadrilateral!trirectangular} figure on the left (hyperbolic geometry), the fourth angle is acute; in the middle figure it is right, and in the right figure it is obtuse (spherical geometry).

 In the same part, Khayy\=am\index{Khayy\=am@al-Khayy\=am, `Umar} makes a detailed analysis of the consequences of the hypothesis on the other two angles, depending on whether they are right, acute, or obtuse. He proves in particular that the hypothesis of the acute angle implies the existence of two straight lines perpendicular to the same straight line and diverging on both sides of this straight line. This is a well-known phenomenon in hyperbolic geometry (but which Khayyam dismisses, considering that it contradicts a first principle). Gerolamo Saccheri,\index{Saccheri, Giovanni Girolamo} whom we have already mentioned,\footnote{See Footnote \ref{n:Saccheri}.} studied the same quadrilaterals several centuries later, and he obtained the same conclusions.

  Finally, let us mention, among the group of\index{Euclid} Arab mathematicians, Na\d{s}\=\i r al-D\=\i n al-\d{T}\=us\=\i , born in \d{T}\=us\index{Tusi@al-\d{T}\=us\=\i , Na\d{s}\=\i r al-D\=\i n} (in present-day Iran) in 1201 and who died in Baghdad in 1274,\footnote{Youschkevitch (\cite{Youschkevitch} p. 414) considers that the \emph{Treatise on the complete quadrilateral} of Al-\d{T}\=us\=\i \ made its author ``the most illustrious scholar of Islam in the field of trigonometry".} of which we have already spoke. Two of his treatises on the theory of parallels are known, the \emph{Treatise which delivers from doubts concerning parallel lines} and his \emph{Redaction of Euclid}, of which a Latin translation and commentary are contained in Vol. 2 of the \emph{Complete Works} of John Wallis\index{Wallis, John}, published in Oxford in 1693. Like `Umar al-Khayy\=am\index{Khayy\=am@al-Khayy\=am, `Umar}, Na\d{s}\=\i r al-D\=\i n al-\d{T}\=us\=\i \ made an\index{Tusi@al-\d{T}\=us\=\i , Na\d{s}\=\i r al-D\=\i n}  analysis of the birectangular isosceles quadrilaterals\index{quadrilateral!birectangular isosceles} which we already discussed.\index{birectangular isosceles quadrilateral} These quadrilaterals will be mentioned again later in this chapter.

   \subsection{The Renaissance}
\begin{figure}[htbp]
\centering
\includegraphics[width=1\linewidth]{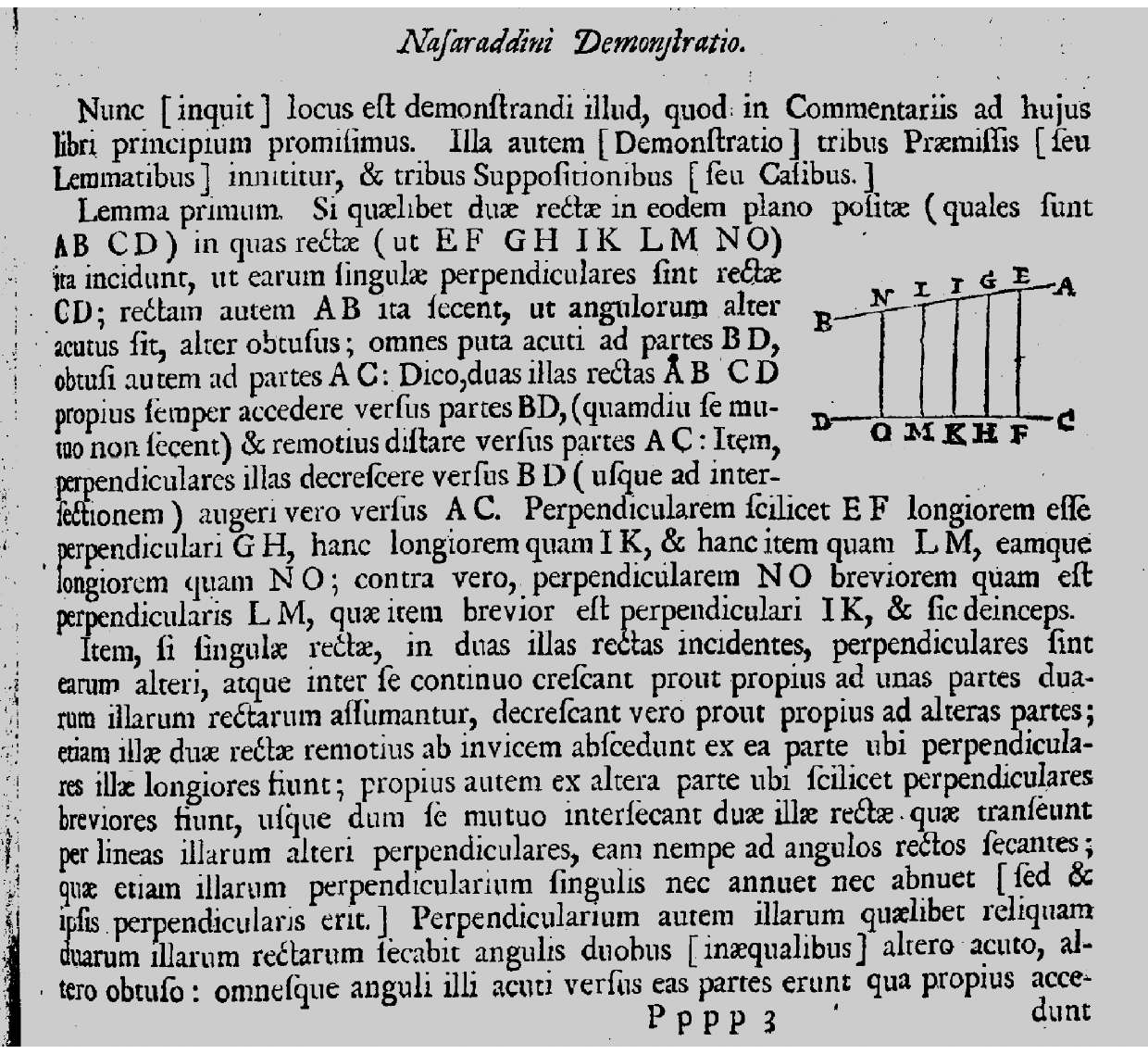}
\caption{\small Page extracted from Wallis' \emph{Complete Works}, a section titled \emph{Nasaraddini Demonstratio} in which the latter sets out the proof of Na\d{s}\=\i r al-D\=\i n al-\d{T}\=us\=\i .}
\label{fig:Wallis-Nasir_gris}
\end{figure}

  At the Renaissance, several mathematicians in Western Europe were naturally interested in the problem of parallels. Among them we mention Christophorus Clavius (1538-1612), Giovanni Alfonso Borelli (1608-1679), John Wallis (1616-1703) and Vitale Giordano da Bitonto (1633-1711). 
  All these mathematicians were familiar with the Greek writings on the subject (especially the work of Proclus), and some were also familiar with the writings of the Arabs. The writings of Na\d{s}\=\i r al-D\=\i n al-\d{T}\=us\=\i \ were
known in Europe thanks to their Latin translation published by Wallis (see Figure \ref{fig:Wallis-Nasir_gris}). Let us also recall that this translation contains a summary of the work of `Umar al-Khayy\=am\index{Khayy\=am@al-Khayy\=am, `Umar} 
and Ibn al-Haytham 
on parallels. 
It is also known that Borelli\index{Borelli, Giovanni Alfonso} collaborated with Abraham Ecchelensis\index{Ecchelensis, Abraham} (1605-1664), a Christian Arab (Syro-Lebanese) Catholic philosopher and linguist, to translate from Arabic manuscripts a commentary on books V and VII of Apollonius' \emph{Conics}.\index{Apollonius of Perga}

Borelli had access to a large Arabic mathematical literature. The proof he gave of the parallel postulate is essentially that of Ibn Qurra and Ibn al-Haytham, see  Rosenfeld, \cite[p. 96]{Rosenfeld}.

The German Jesuit mathematician and physicist Christophorus Clavius,\index{Clavius, Christophorus} translated Euclid\index{Euclid} into Latin and attempted at least two proofs of the fifth postulate, one of which is based on an argument similar to that of Na\d{s}\=\i r al-D\=\i n al-\d{T}\=us\=\i . Clavius gave this proof in 1574 under the assumption that a curve coplanar to a straight line and equidistant from it is also a straight line. Clavius' proof is analyzed by Pont in \cite[p. 195-199]{Pont}. Borelli \cite{Borelli} and Giordano da Bitonto\index{Giordano da Bitonto, Vitale} \cite{Bitonto} showed that Euclid's axioms, without the one on parallels, imply this equidistance property, after introducing a hypothesis that is equivalent to the parallel axiom, see \cite{Bonola-2} p. 13 and \cite{Pont} p. 200-203 and 368. Giordano also proved that if three points on a given line are equidistant from another line, then the two lines are equidistant.

Wallis,\index{Wallis, John} for his part, showed that the parallel postulate\index{parallel postulate} is equivalent to the assertion that there exist, for a given triangle, homothetic triangles that are not congruent to it and that have an arbitrary size \cite{Wallis22}. Having considered, like his predecessors, that the parallel postulate deserved a thorough examination, Wallis proposed that, in the teaching of geometry, this postulate be replaced by a postulate stipulating the existence of homothetic figures for arbitrary figures (and not only for triangles).\footnote{\emph{Opera}, vol. II, Oxford, 1693, p. 667.} Wallis justified\index{Euclid!\emph{Elements}} this idea by the fact that one of Euclid's postulates states that there exist homothetic\index{Euclid} circles of any size.\footnote{Third postulate, Heath's edition of the \emph{Elements} \cite[p. 199]{Heath3}: ``To describe a circle with any centre and distance.'' (A more faithful translation would be: ``And [let it be postulated] that with any centre and any distance a circle can be drawn"; in Greek: ``kai panti kentro kai diast\'emati kuklon gr\'afestai".}

Saccheri, whom we have already mentioned, is considered, along with Lambert and Legendre (whom we will discuss later), as one of the precursors of hyperbolic geometry\index{hyperbolic geometry} who came closest to it. Beltrami,\index{Beltrami, Eugenio} in an article on Saccheri titled {\it Un precursore italiano di Legendre e di Lobatschewski} \cite{Beltrami-Precursore}, notes that several results attributed to Legendre and Lobachevsky were already contained in Saccheri's memoir\index{Saccheri, Giovanni Girolamo} {\it Euclides ab omni n\oe vo vindicatus} (Euclid freed of every flaw). Among these there is the result affirming that in neutral geometry\index{neutral geometry} the sum of the angles in a triangle is smaller than two right angles. This memoir also contains the result saying  that the parallel postulate\index{parallel postulate} can be replaced by the postulate stating the existence of homothetic and non-congruent triangles. Saccheri's memoir was analyzed by several mathematicians after Beltrami, including Mansion \cite{Mansion},
Segre \cite{Segre} and Veronese \cite{Veronese}, see also the recent edition by De Risi \cite{Saccheri}.
   Saccheri is best known for his analysis of quadrilaterals with two adjacent right angles on the same side and the other two angles congruent (Figure \ref{fig:2-Lambert-11}). Saccheri's memoir was published in 1733, the year of its author's death. It is not known whether Saccheri saw his work printed. The memoir was almost forgotten for over 150 years, until it was discovered in 1889 by the Jesuit priest Manganotti, who communicated it to Beltrami.\footnote{The memoir is cited, however, in Kl\"ugel's dissertation (1763) \cite{Klugel}, and Lambert was familiar with this thesis.}

 G. B. Halsted, in the introduction to a translation he published of Saccheri's treatise (see \cite{Saccheri}), reports that Borelli\index{Borelli, Giovanni Alfonso} proposed\index{parallel lines!definition!Saccheri}  the following definition in 1658: ``Two lines are parallel if they are coplanar and have a common perpendicular". Halsted also reports on an attempt of  proof of the parallel postulate by Robert Simson\index{Simson, Robert} in 1756, which starts from the hypothesis that a straight line cannot approach another and then move away from it without intersecting it.\footnote{Such a property of monotonicity is contained in the work of Khayy\=am\index{Khayy\=am@al-Khayy\=am, `Umar} (see \cite{Rashed0} p. 283).}
  
As an example, let us consider Wallis' proof\index{Wallis, John} of the parallel postulate,\index{parallel postulate} which is based on an axiom asserting the existence of homothetic triangles. It is contained in Wallis' article, and is reproduced by St\"ackel and Engel in \cite{Engel-Staeckel}, p. 21-30. It is also reported on by Pont in \cite{Pont} and by Bonola in \cite{Bonola-2}.

\begin{figure}[htbp]
\centering
\includegraphics[width=1\linewidth]{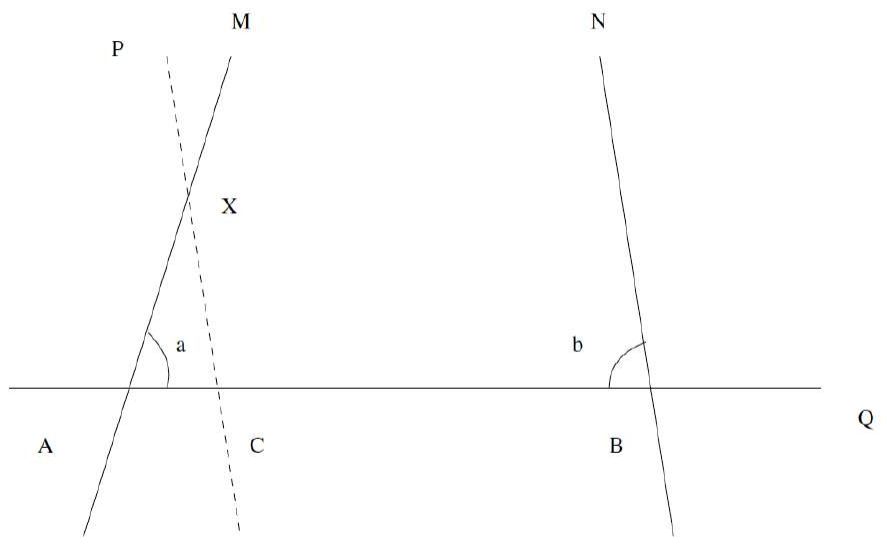}
\caption{\small Figure illustrating Wallis' ``proof" of the parallel postulate.}
\label{fig:wallis1}
\end{figure}

 Consider in the neutral plane\index{neutral plane} (i.e. the plane geometry in which all\index{Euclid} of Euclid's axioms are satisfied except that of parallels, the latter being deleted from the axioms without being replaced by another axiom) two lines $l_1$ and $l_2$, meeting a third line $l_3$ at two points $A$ and $B$ and making with it at these two points, and on the same side, angles $a$ and $b$ satisfying $a+b<\pi$. (These are lines considered in the statement of the parallel axiom; see Figure \ref{fig:wallis1}.)

We move continuously the point $B$ on the line $l_3$ towards the point $A$, by a motion during which the line $l_2$ moves at the same time while keeping a constant angle with $l_3$. When the point $B$ reaches a position $B'$ close enough to $A$, the line $l_2$ (in its new position, which we call $CP$) necessarily meets the line $l_1$. (This is a result of neutral geometry,\index{neutral geometry} which Wallis established before this proof.) Let $C$ be the intersection point of $l_4$ with $l_1$. We consider the triangle $AB'C'$ and we follow the construction that we have just made, in reverse order. That is, now we move the point $B'$ continuously towards $B$, the line $l_4$ being carried by the point $B'$ during this motion in such a way that it always makes the same angle with line $l_3$. The axiom affirming the existence of homothetic triangles says that there exists a triangle $ABC$ homothetic to $AB'C'$ whose
side $AB'$ carries the line $l_3$, and such that the triangle $ABC$ is located on the same side as $AB'C$ with respect to line $l_3$. This shows that lines $l_3$ and $l_2$ intersect (at the point $C$).
   
   This argument is problematic and in fact, in hyperbolic\index{hyperbolic geometry} geometry (which is part of neutral geometry),\index{neutral geometry} when we move the point $B'$ and the line $l_4$ that contains it so that it always maintains the same angle with the line $l_3$, this line ends up no longer meeting the line $l_1$.

 \subsection{The modern period} 
 The modern period in mathematics can be considered to begin with Leonhard Euler. We have mentioned in \cite{NP-Solid-Angles} two attempts of the latter, as well as one by Lagrange,  to prove the parallel postulate.\index{parallel postulate} Let us also quote d'Alembert,\index{Alembert@d'Alembert, Jean le Rond} who, in his works on mathematics, called for simple definitions that are within everyone's reach. He declared that the problem with the question of parallel lines lies in the definitions, writing his famous sentence, that ``the definition and properties of the straight line are the pitfall and the scandal of Geometry", see  \cite[p. 318]{Alembert-Essai}.

Let us mention now an attempt by  Abraham Gotthelf K\"astner\index{K\"astner, Abraham Gotthelf} (1719-1800). K\"astner was appointed professor of mathematics at the University of G\"ottingen in 1756, the year Lambert arrived in that city for a stay that was to last two years. The two men met and remained in contact by letter until Lambert's death. K\"astner was interested in the question of parallels, and he himself attempted to prove the fifth postulate. It is likely that K\"astner contributed to awakening Lambert's interest in the subject. He had a reputation as a very good teacher and was known for his introductory lectures on various subjects. He supervised the doctoral thesis of Georg Simon Kl\"ugel,\index{Kl\"ugel, Georg Simon} defended in 1763 \cite{Klugel}, in which the latter analyzes twenty-eight attempts to prove the parallel axiom made by his predecessors, including that of Saccheri\index{Saccheri, Giovanni Girolamo}, exposing the error contained in each of them. In conclusion of his thesis, Kl\"ugel puts forward the hypothesis that the postulate is not provable. He writes in particular (translation of Pont \cite{Pont} p. 463 et seq.):

 \begin{quote}\small
 Some have tried to exclude this proposition from the list of axioms and have tried to prove it. But these attempts are not free from errors. Other axioms, neither simpler nor more certain, have often been used instead of Euclid's.\index{Euclid} It thus appears, when we consider all these attempts, that Euclid did the right thing by including among the axioms a proposition that could not be properly established with the help of the other propositions.
\end{quote}
 Lambert held the same opinion regarding his predecessors; see \S 3 of his \emph{Theorie der Parallellinien}. He was familiar with Kl\"ugel's thesis,\index{Kl\"ugel, Georg Simon} which he quotes in his memoir.

 Gauss,\index{Gauss, Carl Friedrich} during his studies in G\"ottingen, attended K\"astner's lectures (but not assiduously; it is said that he found them too elementary). Other students of K\"astner included J. M. C. Bartels, who later\index{hyperbolic geometry} became a professor at Kazan University where his lectures were attended by Lobachevsky,\index{Lobachevsky, Nikolai Ivanovich} and Wolfgang Bolyai\index{Bolyai, Wolfgang} (1775-1856),\footnote{Wolfgang Bolyai is the father of J\'anos Bolyai\index{Bolyai, J\'anos}, one of the three inventors of hyperbolic geometry. He was himself a   talented mathematician, poet, and violinist, and a close friend of Gauss.\index{Gauss, Carl Friedrich} He was interested in the foundations of geometry and, like many geometers before him, he tried to prove the parallel postulate. He eventually concluded that this was an extremely difficult problem, and when he learned that his son Janos was working on it, his first reaction was to advise him to abandon this research, see \cite{Meschkowski, Gray, Gray-idea, Bonola-2, Greenberg}. He is credited with a proof that the parallel axiom is equivalent to the statement that three non-aligned points always lie on a circle.
 Rosenfeld points out in \cite[p. 109]{Rosenfeld} that Wolfgang Bolyai, who had not realized the importance of his son's discovery, himself published, after this discovery (and therefore after Lobachevsky's), a proof of the parallel postulate.} who also attended his lectures\index{K\"astner, Abraham Gotthelf} in G\"ottingen. Thus,   K\"astner had a direct or indirect influence on four of the greatest names in the theory of parallels: Lambert and the three founders of non-Euclidean geometry.

 The mathematician Johann Friedrich Lorenz (1738-1807) gave a proof of the equivalence between the parallel axiom and the statement saying that a line passing through a point interior to an angle whose value is strictly between 0 and 180 degrees necessarily meets one of the two sides of the angle \cite{Lorenz}.

 The equivalence between the parallel axiom as stated in the \emph{Elements} and\index{Euclid!\emph{Elements}} the statement that for any line $l$ and for any point $A$ that is not on $l$, there exists a line $l'$ containing $A$ and disjoint from $l$, is attributed to the mathematician John Playfair (1748-1819), even though this axiom was already contained in Ibn al-Haytham's {\it Book on the resolution of doubts raised by Book I of\index{Euclid} Euclid's Elements}\index{Ibn al-Haytham} which we already mentioned, cf. Pont \cite[ p. 171]{Pont}.
      
      The equivalence between the parallel axiom and the fact that there is a quadrilateral with four right angles is contained in the work of Alexis-Claude Clairaut (1713-1765).

Pierre-Simon de Laplace (1749-1827), in his \emph{Exposition of the System of the World}, proposed replacing the parallel axiom with an axiom of similarity, much more natural than Euclid's\index{Euclid}, see \cite[p. 386]{Laplace} and \cite[Vol. VI, p. 472]{Laplace-Oeuvres}. Lazare Carnot\index{Carnot, Lazare} (1753-1823) made the same proposal in his \emph{G\'eom\'etrie de position} (\cite{Carnot} p. 481, in a note).

 Adrien-Marie Legendre\index{Legendre, Adrien-Marie} (1752-1833) is also known for his work on the theory of parallels.
In the third edition  of his \emph{Elements of Geometry} (published in 1800), he provides a proof of the theorem saying that in neutral geometry\index{neutral geometry} the sum of the angles in any triangle is $\leq \pi$, without using the parallel postulate \cite[Proposition 19 of Book I]{Legendre}. This theorem, which is already found in Saccheri's memoir, is often still called Legendre's first theorem.\index{Legendre's first theorem} The reason this result is attributed to Legendre is that, as already noted, Saccheri's work had been forgotten for more than 150 years, and only reappeared in 1889.
Saccheri's proof of this theorem is reproduced in Bonola \cite[p. 56 of the English translation]{Bonola-2}. Gauss also states the Saccheri--Legendre theorem, with a proof analogous to that given by Legendre (see \cite{Gauss}, Vol. VIII). The reader is referred to \cite{ACP} for the place this theorem occupies in the development of neutral geometry and\index{neutral geometry} for an exposition of the proof Legendre gave of it.

  It should be noted here that in the same edition of his \emph{Elements of Geometry}, Legendre introduced a ``proof" of the fact that this sum cannot be strictly less than two right angles, again without using the axiom of parallels. (We know that this proposition is equivalent to the axiom of parallels). Legendre himself acknowledged his mistake later, in an article he published in 1833, the year of his death; see \cite[p. 371-372]{Legendre1833}.  

Legendre\index{Legendre, Adrien-Marie}  is considered, along with Lambert and Saccheri,
\index{Saccheri, Giovanni Girolamo}  to have been one of the closest to the discovery of hyperbolic geometry.\index{hyperbolic geometry} It is interesting to read what Lobachevsky\index{Lobachevsky, Nikolai Ivanovich}   himself wrote about Legendre  \cite[p. 7]{Loba-Nouveaux}:

\begin{quote}\small
 
In the 1833 \emph{M\'emoires de l'Acad\'emie de France}, [Legendre] adds the proposition that the sum of the angles is equal to $\pi$ in all triangles if it has this value in a single one. I also had to prove the same thing in my theory written in 1826. I even find that Legendre sometimes stumbled upon the path I so happily followed; but no doubt a prejudice in favor of the universally accepted proposition led him, at every step, to conclude hastily or to introduce prematurely notions that should no longer be admitted in the new hypothesis.
\end{quote}

 The proposition about the sum of the angles of triangles that Lobachevsky talks about is already found in Saccheri's memoir.

Another prominent mathematician who attempted to prove the parallel postulate is  Joseph Fourier\index{Fourier, Joseph} (1768-1830). Between 1820 and 1827, he wrote a 250-page manuscript containing nearly twenty attempts to prove this postulate. Fourier's manuscript is analyzed\index{parallel lines!definition!Fourier}  at length by Pont in \cite[p. 531-586]{Pont}. He used the following definition of parallelism:
\begin{quote}\small
A line $l$ is parallel to another line $l'$ located in the same plane if $l$ separates the lines in the plane that intersect $l'$ from those that do not.
\end{quote}
This is interesting because this definition is very close to the one used by\index{Lobachevsky, Nikolai Ivanovich}  Lobachevsky, Bolyai\index{Bolyai, J\'anos} and Gauss\index{Gauss, Carl Friedrich} in the setting of hyperbolic geometry.\index{hyperbolic geometry} Fourier then\index{parallel postulate} attempted to  prove the parallel postulate by using the fact that a set of points equidistant from a line is the union of two lines (cf. \cite{Pont} p. 554).
 
 Let us conclude this subsection on the theory of parallels in the modern period by quoting Arthur Schopenhauer\index{Schopenhauer, Arthur} (1788-1860) who was of the opinion that the mathematicians' attempts to prove the parallel postulate was pointless. He writes \cite[Vol. 2 Chap. 13]{Schopenhauer}:
   \begin{quote}\small
  The Euclidean method of proof has brought forth from its own womb its most striking parody and caricature in the famous controversy over the theory of parallels, and in the attempts, repeated every year, to prove the eleventh axiom (also known as the fifth postulate). The axiom asserts indeed through the indirect criterion of a third intersecting line, that two lines inclined on each other (for this is the precise meaning of ``less than two right angles"), if produced far enough, must meet. Now this truth is supposed to be too complicated to pass as self-evident, and therefore needs a proof; but no such proof can be produced, just because there is nothing more immediate. 
  \end{quote}

\section{Excursus 5: On quadrilaterals in the theory of parallels}\label{excursus:quadrilaterals}

We have already mentioned that attempts to prove the parallel postulate
from the notion of equidistance gave rise to the analysis of two remarkable types of quadrilaterals, namely, trirectangular quadrilaterals,\index{trirectangular quadrilateral}\index{Ibn al-Haytham--Lambert quadrilateral} also called Ibn al-Haytham--Lambert quadrilaterals, and\index{birectangular isosceles quadrilateral} birectangular isosceles quadrilaterals (or the Khayy\=am--Saccheri quadrilaterals).\index{Khayy\=am--Saccheri quadrilateral} Because\index{quadrilateral!birectangular isosceles} of\index{quadrilateral!trirectangular} their importance, we summarize here the history and some properties of these two types of quadrilaterals. A detailed analysis which shows their place in the development of non-Euclidean geometry is made in \cite{ACP}.

 \begin{figure}[htbp]
\centering
\includegraphics[width=1\linewidth]{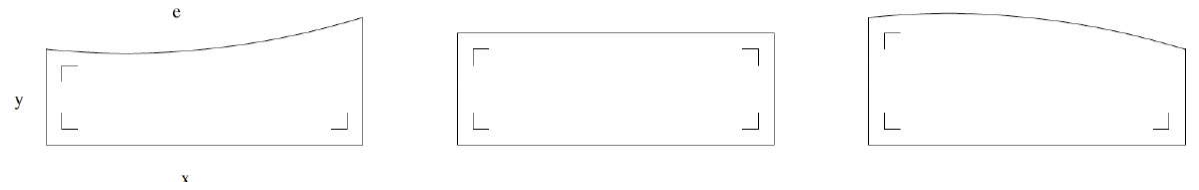}
\caption{\small {Figures (i), (ii), (iii) represent a trirectangular quadrilateral in hyperbolic, Euclidean and spherical geometry respectively .}}
\label{fig:2-Lambert10}
\end{figure}

A {\it trirectangular quadrilateral} in the neutral\index{trirectangular quadrilateral} plane\index{quadrilateral!trirectangular} (hyperbolic or Euclidean) or\index{hyperbolic geometry} on the sphere, is a quadrilateral with three right angles. We have already encountered such a quadrilateral\index{Euclid!\emph{Elements}!Arabic manuscripts} in the works of Ibn al-Haytham, namely, in his {\it Book Explaining the Postulates of Euclid's Book on the Elements}, based on an idea contained in the works of Th\=abit Ibn Qurra\index{Ibn Qurra, Th\=abit}. We also saw that it appeared later in Lambert's  \emph{Theorie der Parallellinien}. These authors made a thorough analysis of the geometry of these quadrilaterals, distinguishing three cases for the remaining angle: acute, right or obtuse (see Figure \ref{fig:2-Lambert10}). Ibn al-Haytham and Lambert both tried to show that the existence of such quadrilaterals where the fourth angle is acute or obtuse contradicts the axioms of neutral geometry.\index{neutral geometry} In the case of the obtuse angle, it is known that this indeed contradicts the axioms of this geometry --- this is Legendre's theorem already mentioned --- and the conclusion and arguments of these authors are correct. Trirectangular quadrilaterals exist in spherical geometry where other axioms of Euclid are not satisfied. In the case of the acute angle, Ibn Qurra and Lambert tried to find a property contradicting the axioms of neutral geometry,\index{neutral geometry} but they could not succeed,\index{hyperbolic geometry} since such quadrilaterals do indeed exist in hyperbolic geometry.

   The history of birectangular isosceles quadrilaterals is similar.

\begin{figure}[htbp]
\centering
\includegraphics[width=1\linewidth]{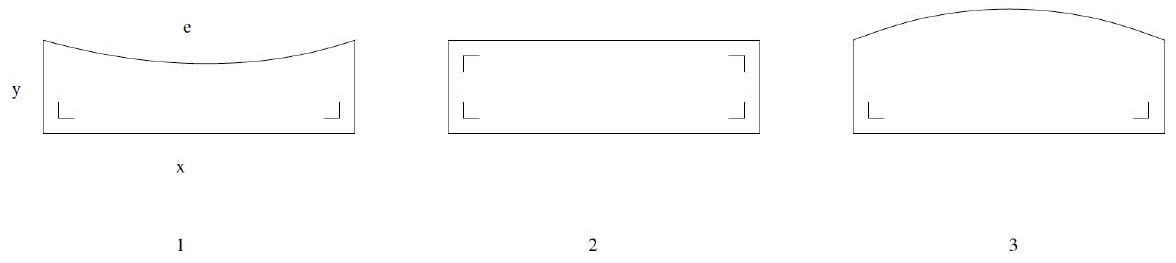}
\caption{\small {Figures (i), (ii), (iii) represent an birectangular isosceles quadrilateral in hyperbolic, Euclidean and spherical geometry respectively.}}
\label{fig:2-Lambert12}
\end{figure}

A birectangular isosceles quadrilateral\index{birectangular isosceles quadrilateral} (also called a Khayy\=am--Saccheri quadrilateral), in\index{quadrilateral!birectangular isosceles} the neutral plane or on the sphere, is a quadrilateral with two opposite sides congruent and forming a right angle with a third side (Figure \ref{fig:2-Lambert12}). In the Euclidean plane, such a\index{quadrilateral!trirectangular} quadrilateral is (as is the case for trirectangular quadrilaterals) a\index{trirectangular quadrilateral} rectangle. By cutting a birectangular isosceles quadrilateral along an axis of symmetry, we obtain two trirectangular quadrilaterals.\index{trirectangular quadrilateral}\index{Ibn al-Haytham--Lambert quadrilateral}

\begin{figure}[htbp]
\centering
\includegraphics[width=0.6\linewidth]{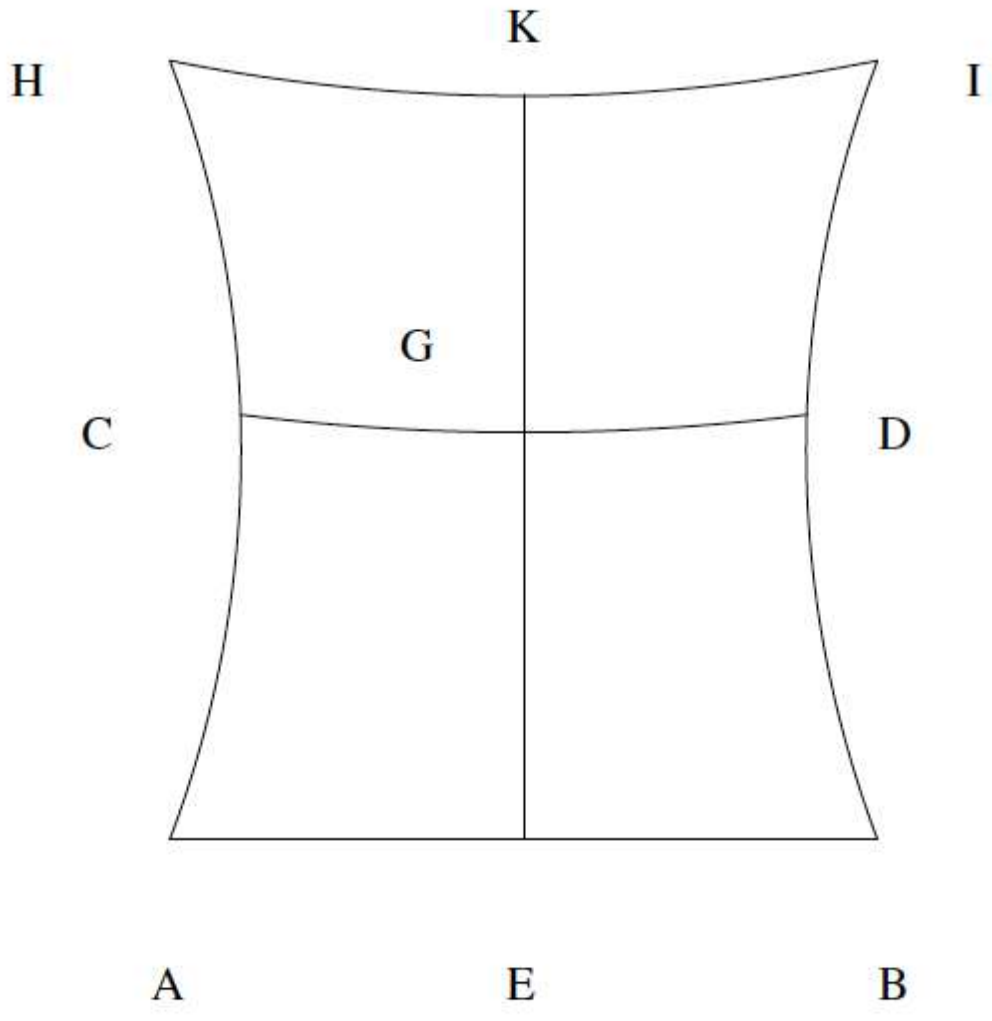}
\caption{\small In his analysis of the Khayy\=am--Saccheri quadrilateral (with right angles at $A,B$) Khayy\=am considers that $AC$ and $BD$ necessarily intersect the perpendicular  at the point $K$ to $EG$, which is not true.}
\label{fig:2-Lambert-11}
\end{figure}

   We have already mentioned that birectangular\index{Khayy\=am--Saccheri quadrilateral} isosceles quadrilaterals were\index{quadrilateral!birectangular isosceles} studied\index{birectangular isosceles quadrilateral} by `Umar al-Khayy\=am\index{Khayy\=am@al-Khayy\=am, `Umar} in Book I of his {\it Commentary on the Difficulties of Some Postulates in Euclid's Work}, cf. \cite{Rashed0} p. 306\index{Euclid} et seq. and \cite{Rosenfeld-Yousch} p. 467. They were also studied by Na\d{s}\=\i r al-D\=\i n al-\d{T}\=us\=\i \ in his two treatises, \emph{The treatise that delivers from the doubt concerning parallel lines} and the \emph{Redaction of Euclid}, cf. \cite{Heath3} Volume I, p. 209 and \cite{Rosenfeld-Yousch} p. 468. They were rediscovered by Saccheri\index{Saccheri, Giovanni Girolamo}, in his {\it Euclides ab omni n\oe vo vindicatus} \cite{Saccheri} (1733) who, like his predecessors, analyzed them by distinguishing three cases, depending on whether the remaining angles (which are equal to each other) are acute, right or obtuse. As with the case of \index{trirectangular quadrilateral} trirectangular quadrilaterals, all these authors considered that the cases of acute\index{quadrilateral!trirectangular} and obtuse angles are excluded from neutral geometry.\index{neutral geometry} Lambert, in his \emph{Theorie der Parallellinien}, reasons several times on trirectangular quadrilaterals obtained by cutting them into two isosceles\index{birectangular isosceles quadrilateral} birectangular quadrilaterals,\index{quadrilateral!birectangular isosceles} notably when he studies the properties of monotonicity of the lengths of the perpendiculars drawn from the line containing the side adjacent to the angle which is not assumed to be a right angle a priori, and those of the angles which these perpendicular segments make with the line from which they are drawn.  We have already discussed
Lambert's analysis. Let us say a few words on the analysis of `Umar al-Khayy\=am\index{Khayy\=am@al-Khayy\=am, `Umar}. It is made in Proposition 3 of Book I of his {\it Commentary} (p. 324 of the translation in \cite{Rashed0}). Khayy\=am considers the quadrilateral $ABCD$ where the angles at $A$ and $B$ are right. The proposition says that the angles $C$ and $D$ are also right angles. (In Proposition 1 of the same book, Khayy\=am had already proved that the angles at $C$ and $D$ are equal.) The proof of Proposition 3 begins with a construction that is erroneous. Indeed, Khayy\=am takes the perpendicular bisector of the segment $AB$ from $E$, the midpoint of $AB$ (Figure \ref{fig:2-Lambert-11}). This perpendicular bisector cuts the segment $CD$ at a point he calls $G$, and he extends it to a point $K$ satisfyings $GK=EG$. Through the point $K$, he draws the perpendicular $HKI$ to the line $KE$. He then extends the lines $AC$ and $BD$ and says that they meet $HKI$ at two points he calls $H$ and $I$. This is not correct, because in neutral geometry\index{neutral geometry} (and in particular in hyperbolic geometry) the\index{hyperbolic geometry} lines $AC$ and $BD$ do not necessarily intersect the line $HKI$.\footnote{Khayy\=am\index{Khayy\=am@al-Khayy\=am, `Umar}, in his ``proof" (which is defective), passes over a number of details, with statements such as: ``This is one of the things that have already been treated by the Philosopher" (i.e., Aristotle\index{Aristotle}); ``This is one of the things that you will be able to recognize with a minimum of reflection and investigation", or again, ``It will be recognized with a minimum of reflection and we will therefore omit it in order to avoid prolixity. Thus, whoever wishes to establish this here according to the mathematical order, let him do so; we will not prevent him!" (quoted by Rashed and Vahabzadeh in \cite{Rashed0} p. 326).}

    Saccheri, \index{Saccheri, Giovanni Girolamo} in the memoir in \cite{Bonola-2} we already mentioned, analyzes the three possibilities for the angles at $C$ and $D$ (which he knows to be congruent): acute, right, or obtuse. He proves that if one of the three hypotheses holds for\index{birectangular isosceles quadrilateral} a particular\index{quadrilateral!birectangular isosceles} birectangular isosceles  quadrilateral,\index{Khayy\=am--Saccheri quadrilateral} then it holds for all birectangular isosceles quadrilaterals (Propositions V, VI, and VII). The obtuse angle hypothesis is ruled out by an argument using the work of Na\d{s}\=\i r al-D\=\i n al-\d{T}\=us\=\i ,\index{Tusi@al-\d{T}\=us\=\i , Na\d{s}\=\i r al-D\=\i n} which is quoted by Saccheri.

Saccheri\index{Saccheri, Giovanni Girolamo} then proves that the right, obtuse, and acute angle hypotheses imply that the sum of the angles in a triangle is equal to $\pi$, greater than $\pi$, and smaller than $\pi$, respectively. In doing so, he establishes a series of results valid in neutral geometry,\index{neutral geometry} with the aim of finally proving that the acute angle hypothesis cannot hold. To conclude his argument, and having not reached a proof, he notes that the acute angle hypothesis implies the existence of perpendiculars ``located at infinity" , that is, the existence of two disjoint lines having no common perpendicular, a property that he excludes \emph{a priori}. We know from the later work of Lobachevsky\index{Lobachevsky, Nikolai Ivanovich} that this situation arises in hyperbolic geometry.\index{hyperbolic geometry}

 \bigskip
 
Author's addresses

 Athanase Papadopoulos (corresponding author),  Institut de Recherche Math\'{e}matique Avanc\'{e}e,   
 Universit\'{e} de Strasbourg et CNRS,   
  7 rue Ren\'{e} Descartes,  67084 Strasbourg Cedex, France.
 email: athanase.papadopoulos@math.unistra.fr

 \bigskip

 Guillaume Th\'eret,   Laboratoire Interdisciplinaire Carnot de Bourgogne (ICB),
9 avenue Alain Savary,
BP 47870, 21078 Dijon Cedex
FRANCE.
email: guillaume.theret71@orange.fr

\printindex

\begin{thebibliography}{99}


\bibitem{ACP} N. A'Campo and A. Papadopoulos, Notes on Hyperbolic Geometry. In: \emph{Strasbourg Master class on Geometry} (ed. A. Papadopoulos),  IRMA Lectures in Mathematics and Theoretical Physics, vol. 18, Z\"urich,  European Mathematical Society (EMS) p.~1-182, 2012.
   
  \bibitem{Alembert-Essai} d'Alembert (J. Le Rond d'Alembert),
    Essai sur les \'el\'ements de philosophie, first edition, 1759, 
    new edition,  Fayard, Paris (1986)
 
\bibitem{Barbarin1928} P. Barbarin, \textit{La G\'eom\'etrie non-Euclidienne
 (suivie de Notes Sur la G\'eom\'etrie non-Euclidienne dans ses Rapports avec la Physique Math\'ema\-tique par A. Buhl)}. Gauthier-Villars, Paris,  3rd ed. 1928. (1st ed. 1902.)

 
\bibitem{Barbarin-Correspondance} P. Barbarin, La correspondance entre Ho\"uel et de Tilly,  \textit{Bulletin des Sciences Math\'ematiques}. 50 (1926), p. 50-61.
  
  \bibitem{Beltrami-Precursore} E. Beltrami, Un precursore italiano di Legendre e di Lobatschewski. \textit{Rendiconti della Reale Accademia dei Lincei  Roma} (4)    5 (1889),  p.~441-448. \emph{\OE vres compl\`etes de Beltrami},  Vol. IV p.~348.
 
\bibitem{Belyi} Y. A.  Belyi, {\it Sur un ouvrage d'Euler de g\'eom\'etrie} (Russian), 14 (1961), p.~237-284.
   
\bibitem{Biot} J.-B. Biot, \emph{M\'elanges scientifiques et litt\'eraires}, t. II, Michel  L\'evy Fr\`eres, Libraires-\'editeurs, Paris,  1858.

\bibitem{Beltrami-Boi} L. Boi, L. Giacardi and R. Tazzioli (ed.), {\it La d\'ecouverte de la g\'eom\'etrie non-euclidienne sur la pseudosph\`ere:  Les lettres d'Eugenio Beltrami \`a Jules Ho\"uel},
        Introduction, notes et commentaires par L. Boi, L. Giacardi et R. Tazzioli.
Preface by Ch. Houzel and E.~Knobloch, Sciences dans l'Histoire,  Librairie Scientifique et Technique Albert Blanchard, Paris, 1998.
 
\bibitem{Bolyai} J. Bolyai, Scientiam Spatii absolute veram exhibens; a veritate aut falsitate axiomatis XI Euclidei (a priori haud unquam decidenda) independentem; adjecta ad casum falsitatis quadratura circuli geometrica, 1832. 
Appendix to the \emph{Tentamen} by F. Bolyai. Budapest, 1902.  French translation by J.~Ho\"uel, La Science absolue de l'espace,
ind\'ependante de la v\'erit\'e ou de la fausset\'e de l'Axiome XI, pr\'ec\'ed\'e d'une notice sur la vie et
les travaux de W. et de J.\ Bolyai par M. Fr. Schmidt. \emph{M\'emoires de la Soci\'et\'e des Sciences Physiques et Naturelles de Bordeaux} 5  (1867), p.~189-248, et   Gauthier-Villars, Paris, 1868. 

\bibitem{Bonola-2} R. Bonola,
 La geometria non-euclidea. Esposizione storico-critica del suo sviluppo.  1st ed., Ditta Nicola Zanichelli editore, Bologna, 1906. German translation by M. Liebmann in the collection Wissenschaft und Hypothese, Teubner, Leipzig, 1908. English translation by  H. S. Carslaw,  Non-Euclidean geometry, A critical and historical study of its development. 1st \'ed., Chicago, 1912, reprint  Dover, New York, 1955.
     
     

 \bibitem{Borelli} G. A. Borelli, \emph{Euclides restitutus, siue, Prisca geometriae elementa},  Pisa, 1658.

\bibitem{Brunel} G. Brunel, Notice sur l'influence scientifique de Guillaume-Jules Ho\"uel,  M\'emoires de la Soci\'et\'e de Sciences Physiques et Naturelles de Bordeaux, 4 (1888), p. 1-78.


\bibitem{Carnot} L. N. M. Carnot, G\'eom\'etrie de position, J. B. M. Duprat, Paris, 1803.



\bibitem{Chatelet} A. Ch\^atelet, Rapports sur la g\'eom\'etrie hyperbolique,  Bulletin des Sciences Math\'ematiques, 2e s\'erie, 37 (1913),  p.  134-144.
 
 
\bibitem{Clavius} C. Clavius, Euclidis Elementorum libri XV, Rome 1574.
    

\bibitem{Morgan} A. De Morgan, A budget of paradoxes, Longmans, Green, London,  1872. Reprint: The Open Court Publishing Co., 1915, et Dover, 1954.


 \bibitem{Engel} F.~Engel, \emph{Nikolaj Iwanowitsch Lobatschefskij}. Zwei geometrische Abhandlungen aus dem russischen \"Ubersetzt mit Anmerkungen und mit einer Biographie des Verfassers. German translation of the  \'El\'ements de g\'eom\'etrie and of the Nouveaux \'el\'ements de g\'eom\'etrie, avec une th\'eorie compl\`ete des parall\`eles.   Teubner, Leipzig, 1898.

%
%
  
\bibitem{Gauss}  C. F. Gauss, Werke, K\"onigliche Gesellschaft der Wissenschaften, G\"ottingen, 1900.


\bibitem{Bitonto} V. Giordano  da Bitonto, Euclide restituito overo gli antichi elementi geometrici ristaurati e facilitati, Rome, 1680.
 


\bibitem{Gray} J. Gray, Worlds out of nothing:
A course in the history of geometry in the 19th century.
 Springer Undergraduate Mathematics Series, Springer-Verlag, London. 1st ed., 2007.
2nd ed., 2010.


\bibitem{Gray-idea} J. Gray, Ideas of space: Euclidean, non-Euclidean and relativistic, Clarendon Press, Oxford, 1979.


\bibitem{Greenberg} M. J. Greenberg, Euclidean and non-Euclidean geometries. 4th ed.,  W. H. Freeman, New York, N.Y., 2008.  (1st ed. 1980).

 \bibitem{Halsted-Monthly} G.B. Halsted, John Henry Lambert (Biography), \emph{The American Mathematical Monthly}, Vol. II, No. 7-8 (1895), p. 208-211.


\bibitem{Heath3} T. Heath, The Thirteen Books of Euclid's Elements translated from the text of Heiberg with introduction and commentary, 3 volumes, First edition 1908; Cambridge University Press, Cambridge, several later editions.

 \bibitem{Heath-Aristotle}  T. Heath,
 Mathematics in Aristotle. Oxford University Press, London, 1949; reprint, Classic Studies in the History of Ideas, Thoemmes Press, Bristol, 1998.
 
 \bibitem{Heiberg} Euclides, Elementa post I. L. Heiberg edidit E. S. Stamatis. (T. I to V of the
Euclidis op\'era omnia.)
Vol. I : Libri I-IV cum appendicibus, vol. II : Libri V-IX cum
appendice. Leipzig, B. G. Teubner, 1883-1888, New printing,  1969-1970, 2 vol. in-8$^{\mathrm{o}}$,
XLii-190 p. et vm-239 p.

\bibitem{Houzel} C. Houzel, The birth of non-Euclidean geometry. In :  1830-1930: A century of geometry,  p. 3-21, Lecture Notes in Physics, Springer, Berlin-Heidelberg, 1992.


\bibitem{Rashed} Ibn al-Haytham,  Les Connus. In: Les math\'ematiques infinit\'esimales du IX${}^e$ au XI${}^e$ si\`ecle, vol. III: Ibn al-Haytham. Th\'eorie des coniques, constructions g\'eom\'etriques et g\'eom\'etrie pratique, London, al-Furq\=an, 446 p.



\bibitem{Klugel} G. S. Kl\"ugel, Conatuum praecipuorum theoriam parallelarum demonstrandi recensio,  quam publico examini submittent Abraham  Gotthelf K\"astner et auctor Georgius Simon Kl\"ugel, G\"ottingen, 1763, 30 p.

\bibitem{Krause} M. Krause, Die Sph\"arik von Menelaos aus Alexandrien, Abhandlungen der Gesellschaft der Wissenschaften zu G\"ottingen, phil.-hist. Klasse, 3, 17, Berlin, 1936.

\bibitem{Lambert-Observations} J. H. Lambert, Observations trigonom\'etriques, M\'em. Acad. sc. Berlin, 24 (1770), p. 327-357.

\bibitem{Lambert-Theorie} J. H. Lambert, Theorie der Parallellinien, manuscript, 1766, published in \cite{Engel-Staeckel} in 1895. French translation in \cite{2012-Lambert}, 2012.

\bibitem{Laplace} P.-S. de Laplace, Exposition du syst\`eme du monde, Imprimerie du Cercle-Social, Paris, 1st ed. 1796.

\bibitem{Laplace-Oeuvres} P.-S. de Laplace, \OE uvres compl\`etes, 14 volumes, publi\'ees sous les auspices de l'Acad\'emie des sciences,  Gauthier-Villars, Paris,  1878-1912.

  \bibitem{Laptev-Theory}  B. L. Laptev, The theory of parallels in the first works of  N. I. Lobachevsky (Russian), Istoriko-Matematicheskie Issledovaniya  (1951) 4, p. 203-229.

 
\bibitem{Laptev-Geometriia} B. L. Laptev, Lobachevsky geometry (Russian),   {e\"e} istoria i znachenie, Moscou, 1976.
     

 \bibitem{Legendre} A.-M. Legendre,  \'El\'ements de g\'eom\'etrie, Firmin Didot, Paris,  First edition 1798; several later editions.
    
 
 \bibitem{Legendre1833} A.-M. Legendre, R\'eflexions sur les diff\'erentes mani\`eres de d\'emontrer la th\'eorie des parall\`eles ou le th\'eor\`eme sur la somme des trois angles du triangle, M\'emoires de l'Acad\'emie des sciences de l'Institut de France, 12 (1833), p.  367-410.
 
        \bibitem{L} N. I.  Lobachevsky, Pang\'eom\'etrie ou pr\'ecis de g\'eom\'etrie fond\'ee sur une th\'eorie g\'en\'erale et rigoureuse des parall\`eles.  In: Recueil d'articles \'ecrits par les professeurs de l'Universit\'e de Kazan \`a l'occasion du cinquantenaire de sa fondation, Vol. I, p.  277-340 (1856).
New edition commented and translated into English by A. Papadopoulos, 
 Heritage of European Mathematics, Vol. 4, European Mathematical Society, Z\"urich, 2010.
 
 
 
\bibitem{Loba-Nouveaux}  N. I.  Lobachevsky, New principles of geometry with a complete theory of parallels (Russian), \textit{Uchenye zapiski kazanskogo imperatorskogo  universiteta} 1835, No. 3, p. 3-48;1836, No. 2, p. 3-98 \& No. 3, p. 3-50;1837, No. 1, p. 3-97;1838, No. 1, p. 3-124 \& No. 3, p. 3-65. 
French translation by par  F. Mailleux,  M\'em. Soc. Royale sci. de Li\`ege (3), 2 No. 5, p. 101;3-101 No. 2, p. 32-32 (1899).

\bibitem{Lorenz} J. F. Lorenz, Grundri\ss \ der reinen und angewandten
Mathematik, oder der erste Cursus der gesamten Mathematik, Helmst\"adt, Fleckeisen 1791,


\bibitem{Mansion} P. Mansion, Analyse des recherches du P. Saccheri, S. J., sur le postulatum d'Euclide,  Annales de la Soci\'et\'e scientifique de Bruxelles, 14 (1889-1890) 2\`eme partie, p. 46-59. 
  
\bibitem{Meschkowski} H. Meschkowski,  Non-Euclidean geometry, English translation by  A. Shenitzer, Acad. Press, New York, 1964.


\bibitem{NP-Solid-Angles} S. Negrepontis and A. Papadopoulos, Notes on angles and solid angles, 
in relation with  Euler's memoir 
 De mensura angulorum solidorum,  This volume.


    
 
  \bibitem{2012-Lambert} A. Papadopoulos and G. Th\'eret, La th\'eorie des parall\`eles de Johann Heinrich Lambert,  critical edition with French translation and mathematical and historical commentary, Librairie Scientifique et Technique Albert Blanchard, coll. Sciences dans l'Histoire, Paris, 214 p., 2014.
  
 

 
\bibitem{Picard} \'E. Picard, La science moderne et son \'etat actuel, Flammarion, Paris, 1914.

\bibitem{Poincare1891} H. Poincar\'e, Les g\'eom\'etries non euclidiennes, Revue g\'en\'erale des sciences pures et appliqu\'ees, 2e ann\'ee, No. 23, 15 d\'ecembre 1891, p. 769-774.

\bibitem{Poincare-Hilbert} H. Poincar\'e, Rapport sur les travaux de M. Hilbert (Rapport relatif au IIIe concours du prix Lobatschefski d\'ecern\'e le 14 f\'evrier 1904). Bulletin de la Soci\'et\'e physico-math\'ematique de Kasan, tome 14, 1904.

 \bibitem{Pont} J.-C. Pont,    L'aventure des parall\`eles. Histoire de la g\'eom\'etrie non-euclidienne~: pr\'ecurseurs et attard\'es. \'ed. Peter Lang, 
 Bern,  1986.

 
\bibitem{Proclus} Proclus, Commentaires sur le Livre I des \'el\'ements d'Euclide, French translation by P. Ver Eecke, Collection des Travaux de l'Acad\'emie Internationale des Sciences, No.~1,  Descl\'ee de Brouwer \& Cie, Bruges, 1948. 

\bibitem{Proclus-Morrow} Proclus, A commentary on the first book of Euclid's Elements,   Translated with Introduction and Notes by Glenn R. Morrow, Princeton University Press, Princeton, New Jersey, 1970.



\bibitem{RR} R. Rashed and A. Papadopoulos,  Menelaus'  Spherics:  Early Translation and al-M\=ah\=an\=\i /al-Haraw\={\i}'s version:
Critical edition of Menelaus' Spherics from Arabic manuscripts, with historical and mathematical commentaries, De Gruyter, Series: Scientia Graeco-Arabica,  21,  2017, 890 pages.     
 

 
 \bibitem{Rashed1}  R. Rashed, Les math\'ematiques infinit\'esimales du IX\`eme au XI\`eme si\`ecle.
     5 volumes, Al-Furqan Islamic Heritage Foundation Publication 109, London.
 Vol. I~: Fondateurs et commentateurs~: Ban\=u M\=usa, Ibn Qurr\=a, Ibn S\=\i n\=an, al-Kh\=azin, al-Qu\d{h}i, Ibn al-S\=am, Ibn Hund, 1996~;
   Vol. II~: Ibn al-Haytham, 1993~;
 Vol. III~: Ibn al-Haytham. Th\'eorie des coniques, constructions g\'eom\'etriques et g\'eom\'etrie pratique, 2000~;
 Vol. IV~: Ibn al-Haytham. M\'ethodes g\'eom\'etriques, transformations ponctuelles et philosophie des math\'ematiques, 2002~;
 Vol. V~:  Ibn al-Haytham. Astronomie, g\'eom\'etrie sph\'erique et trigonom\'etrie, 2006.


 \bibitem{Rashed0}  R. Rashed  et  B. Vahabzadeh,  Al-Khayy\=am math\'ematicien.  Collection Sciences dans l'Histoire,  Librairie Scientifique et Technique Albert Blanchard, Paris, 1999.



\bibitem{RH2005} R. Rashed and C. Houzel, Th\=abit Ibn Qurra et la th\'eorie des parall\`eles, Arabic science and philosophy, 15 (2005), p. 9-55.



 
\bibitem{Rosenfeld} B. A. Rosenfeld,  History of non-Euclidean geometry. Transl. A. Shenitzer, Studies in the History of Mathematics and Physical Sciences 12, Springer-Verlag, New York,  1988.

\bibitem{Rosenfeld-Yousch} B. A. Rosenfeld et A. P. Youschkevitch,  Geometry.  In  Encyclopedia of the history of Arabic science,  ed. R. Rashed and R. Morelon, Vol. 2,  Routledge, London, 1996, p. 447-494.


 \bibitem{Saccheri} G. G. Saccheri,  Euclides ab omni n\oe vo 
 vindicatus: sive conatus geometricus quo stabiliuntur prima ipsa universae geometriae principia. Pauli Antoni Montani, Milan, 1733.  German transl. of Book I (l'ouvrage de Saccheri en contient deux) by P. St\"ackel dans  \cite{Engel-Staeckel}. English transl.: Euclid freed of every flaw: geometrical essay in which are established the fundamental principles of universal geometry by  G. B. Halsted, Chicago and London: the Open court publishing Co., 1920. Italian edition with commentaries by V. De Risi,  Pisa, Edizioni della Normale, 2011. 
   
   \bibitem{Schopenhauer} A. Schopenhauer, The world as will and representation,  Vol. 1, Ch.  Janaway (ed.), Transl.  J. Norman,  A. Welchman,   Cambridge University Press, 2010.
   
   
\bibitem{Segre} C. Segre, Congetture intorno all'influenza di Girolamo Saccheri sulla formazione della geometria non-euclidea,  Atti della Reale Accademia delle scienze di Torino, 38 (1902-1903), p. 535-547.

 \bibitem{Engel-Staeckel}  P.~St\"ackel and F.~Engel,  Die Theorie der Parallellinien von Euklid bis auf Gauss, eine Urkundensammlung zur Vorgeschichte der nicht-euklidischen Geometrie, B. G. Teubner, Leipzig, 1895.

%



\bibitem{Ver-Eecke-Theodosius} Theodosius, Les Sph\'eriques de Th\'eodose de Tripoli, \OE uvres traduites pour la premi\`ere
fois du grec en fran\c cais avec une introduction et des notes par P. Ver Eecke, Bruges,
Descl\'ee de Brouwer, 1927; repr. Paris, Librairie Albert Blanchard, 1959.


%
%
%
   
\bibitem{Veronese} G. Veronese, Fondamenti di geometria a pi\`u dimensioni e a pi\`u specie di unit\`a rettilinee esposti in forma elementare: Lezioni per la scuola di magistero in matematica, Padova,  Tipografia del Seminario, 1891

 
   \bibitem{Wallis22} J. Wallis, Demonstratio Postulati Quinti Euclidis, in Operum Mathematicorum, t. II, Oxford, 1693.
 
%
%


 

\bibitem{Youschkevitch} A. P. Youschkevitch,
 Les math\'ematiques arabes  \emph{(}VIII$^\textrm{e}$-XV$^\textrm{e}$ si\`ecles\emph{)}. Translated from the Russian by M. Cazenave et K. Jaouiche,
 Collection L'Histoire des Sciences, Textes et \'{E}tudes, Librairie Philosophique J. Vrin, Paris, 1976.
\end{thebibliography}
\end{document}